\documentclass[12pt]{article}
\usepackage{graphicx} 
\usepackage[
  pdftex,
  plainpages=false,
  pdfpagelabels,
  colorlinks,
  citecolor=magenta,
  linkcolor=blue,
  urlcolor=black,
  filecolor=black,
  bookmarksopen
]{hyperref} 

\usepackage{amssymb, amsfonts, amsmath, amsthm, enumerate, wasysym}
\usepackage[margin=1in]{geometry}
\usepackage{tikz,etoolbox}
\usepackage{mathtools}
\usepackage[toc,page]{appendix} 
\usepackage{oplotsymbl} 
\usepackage{caption}
\usepackage{bm}
\usetikzlibrary{matrix,calc,shapes,decorations}
\pgfdeclarelayer{background}
\pgfdeclarelayer{main}
\pgfdeclarelayer{foreground}
\pgfsetlayers{background,main,foreground}

\usepackage{graphicx}
\newcommand{\icol}[1]
{                 
\left(\begin{smallmatrix}#1\end{smallmatrix}\right)
}    
\newcommand{\C}{\mathbb{C}}

\newcommand{\Z}{\mathbb{Z}}
\newcommand{\R}{\mathbb{R}}
\newcommand{\Q}{\mathbb{Q}}
\renewcommand{\phi}{\varphi}

\DeclareMathOperator*{\conv}{conv}

\renewcommand{\P}{{\mathcal P}}

\newcommand{\va}{{\mathbf a}}
\newcommand{\e}{{\mathbf e}}

\newcommand{\vol}{\operatorname{vol}}

\newcommand{\K}{\mathcal{K}}

\newtheorem{cor}{Corollary}

\newtheorem{lemma}{Lemma}

\newtheorem{question}{Question}
\newtheorem{eg}{Example}
\newenvironment{example}{\begin{eg}\rm}{\end{eg}}

\newtheorem{definition}{Definition}
\newtheorem{observation}{Observation}
\newtheorem{remark}{Remark}
\usepackage{thmtools}
\usepackage{thm-restate}
\usepackage{cleveref}
\declaretheorem[name=Theorem]{thm}

\title{Ehrhart quasi-polynomials via Barnes polynomials and discrete moments of parallelepipeds}

\author{Sinai Robins}

\date{August 14, 2026}
\begin{document}
\maketitle

\begin{abstract}
We give novel and explicit formulas for the Ehrhart quasi-polynomials of any rational  polytope $\P$, 
in terms of Barnes polynomials and discrete moments of half-open parallelepipeds. 
These formulas  hold for all positive dilations of $\P$. There is an interesting appearance of an auxiliary complex $z$-parameter, which seems to allow for more compact formulations. 

We  give similar formulas for discrete moments of $\P$ and its positive dilates, objects closely related to the literature on sums of polynomials over a polytope.    The appearance of the Barnes polynomials and the Barnes numbers allow for explicit computations.  From this work, it is clear that the complexity of computing Ehrhart quasi-polynomials lies mainly in the computation of discrete moments of  half-open integer parallelepipeds.  These discrete moments are in general summed over a  particular lattice flow on a compact torus, introduced below. 

By developing the properties of discrete moments of integer parallelepipeds further, we realize them as geometric, higher-dimensional Dedekind sums, which we call polytope Dedekind sums.  They possess reciprocity laws that extend the classical Dedekind sum reciprocity law.
     
As an appplication, we find a novel vertex formula for the general codimension-$1$ quasi-coefficient 
$c_{d-1}(t)$ of Ehrhart quasi-polynomials, valid for all $t>0$.
 
Another consequence of the main results involves novel and canonical vanishing identities for a rational polytope $\P$, extending the well-known vanishing identities of Brion-Vergne.  As a further consequence, we obtain a differential equation for discrete moments of $\P$, extending the work of Eva Linke. 

For smooth polytopes,  we obtain novel and much simpler formulations of Ehrhart polynomials, discrete moments, and vanishing identities.   From the perspective of Barnes polynomials and Barnes numbers, such identities may be of independent interest.  These formulations show the utility of Barnes polynomials in geometric combinatorics, due to their very rich structure that extends the $1$-dimensional Bernoulli polynomials.  
\end{abstract}

\newpage
\tableofcontents

\medskip

\section{Introduction, and statements of the main results}  
\label{intro}

One of the most natural definitions for a discrete
volume of a rational polytope $\P \subset \R^d$ is its classical integer point enumerator 
 $L_\P(t) := | \Z^d \cap t\P |$,
 where $t\P$ is the $t$-dilate of $\P$.  When $\P$ is an integer (resp. rational) polytope, and $t$ is a positive integer,  it is known by
Ehrhart's original work \cite{Ehrhart1} that $L_\P(t)$ is a polynomial (resp. quasi-polynomial) in the integer-valued parameter $t$. 

Here we give novel formulas for Ehrhart quasi-polynomials of rational  polytopes, in terms of  Barnes polynomials and discrete moments of half-open parallelepipeds.   These formulas hold for all $t>0$, and for any rational polytope.
For clarity of notation, most of our principal formulas are first stated
for simple rational polytopes, and then extended to all rational polytopes. Corollary
\ref{cor:general rational polytopes} and Observation
\ref{observation:extension-principle-general-polytopes}
show that every result whose
proof uses simplicity only through the simpliciality of the vertex tangent
cones (except of course for smooth polytopes) extends to arbitrary rational polytopes by means of disjoint
partially-open simplicial decompositions of each vertex tangent cone.

We find here two main objects that form the complete building blocks for the Ehrhart quasi-polynomials of rational polytopes:  
the collection of  {\bf Barnes polynomials}  $B_k(t, \va)$, and the {\bf discrete moments of vertex parallelepipeds}, 
defined in \eqref{defining sum of discrete moments} below. 
The Barnes polynomials are
 defined by their generating function:
\begin{equation} \label{Barnes}
\frac{ x^d e^{t x} }{      (e^{a_1 x} - 1) \cdots (e^{a_d x} - 1)    } :=
\sum_{k \geq 0}  
B_k(t, \va) \frac{x^k}{k!},
\end{equation}
where $\va := ( a_1, \dots, a_d ) \in \C^d$,
and where initially $t \in \mathbb \R_{>0}$.  
The identity \eqref{Barnes} is valid for all $x\in \C$ inside the domain of convergence of the series.  
It follows easily from \eqref{Barnes} that for each nonnegative integer $k$,  $B_k(t, \va)$ is a polynomial function of $t$. 
It also follows that when $d=1$,  $B_k(t, (1))$ is the classical Bernoulli polynomial.    
Because the Barnes polynomials are a natural multivariable extension of the classical Bernoulli polynomials, these two families of polynomials enjoy some similar properties. We note that some papers refer to the polynomials $B_k(t, \va)$ as \emph{Bernoulli-Barnes} polynomials, while others refer to them as \emph{multiple Bernoulli} polynomials.   In $1901$ E. W. Barnes \cite{EWBarnes} initiated the systematic study of the polynomials $B_k(t, \va)$.  Barnes also related these polynomials to multiple gamma functions, and to some useful zeta functions. 
In recent years, the Barnes polynomials have been studied extensively by the important work of Bayad et al. \cite{BayadBeck2014, BayadKim2013, BayadSimsek2011}.   

This paper uses elementary techniques, and 
our first main result is an explicit formula for the 
Ehrhart quasi-polynomial of any simple integer rational polytope $\P$, in terms of Barnes polynomials and discrete moments of parallelepipeds.  Such discrete moments may be thought of as belonging to spatial statistics over the integer lattice \cite{Cressie}.  We show that similar formulas also hold for any rational polytope, after triangulating its vertex tangent cones. 

Throughout, we write $V$ for the set of vertices of a polytope $\P$, and we define the usual local information at each vertex $v$ of a simple rational polytope $\P$.  We use the standard notation $d$-polytope for a $d$-dimensional polytope. 

We call an integer vector $w\in \R^d$ 
{\bf primitive} if the $gcd$ of all of its coordinates equals $1$. We let
\[
w_1(v), \dots, w_d(v)
\]
be the (unique) {\bf primitive integer edge vectors} that emanate from a vertex $v\in V$.  We use the column convention
\begin{equation}
\label{def:edge matrix}
W_v
:=
\begin{bmatrix}
w_1(v) & \cdots & w_d(v)
\end{bmatrix},
\end{equation}
so that the primitive edge vectors are the columns of $W_v$.  We also write
\begin{equation}
\label{def:edge lattice}
\Lambda_v:=W_v\Z^d,
\end{equation}
the integer sublattice generated by the edge vectors incident with $v$.

\begin{definition}
For each vertex $v$ of $\P$, the \bf{vertex parallelepiped} $\Pi_v$ is defined by
\begin{equation}
\label{def:vertex parallelepiped}
\Pi_v
:=
W_v[0,1)^d
=
\left\{
\sum_{k=1}^d \lambda_k w_k(v)
\;\middle|\;
0\leq \lambda_k<1 \text{ for } k=1,\ldots,d
\right\}.
\end{equation}
\end{definition}
We note that $\Pi _v$ is by definition a {\bf half-open integer parallelepiped} with {\bf one vertex at the origin}.  We also note that the translated vertex parallelepiped $\Pi _v + v$ tiles the vertex tangent cone $\K_v$ by translations with integral linear combinations of the edge vectors 
$w_1(v), \dots, w_d(v)$.
For any $z\in \C^d$, and $x \in \R^d$, we will use the slightly unusual, complexified inner product 
$\langle x, z \rangle:= \sum_{m=1}^d x_m z_m$, noting the absence of complex conjugation in this definition. 
For each vertex $v \in \P$, we  define the local vector of complex linear forms
\begin{equation}  
\label{def:a_v}
\va_v:= W_v^T z= 
\begin{pmatrix}
 \langle  w_1(v), z \rangle\\
\vdots \\
\langle w_d(v), z \rangle
\end{pmatrix} \in \C^d,
\end{equation}
which is a vector whose $k$'th coordinate is the inner product $\langle w_k, z \rangle$.  
The starting point for some of our proofs is Brion's fundamental result
 for the  integer point transform of a rational polytope $\P$.  For a general \emph{bounded set} $Q\subset \R^d$, we define its {\bf integer point transform} by:
\begin{equation}
\label{def: general def of integer point transform}
\sigma_{Q}(z) := \sum_{q \in  Q  \cap \Z^d}  e^{\langle q, z \rangle}.
\end{equation}
The integer point transform $\sigma_{Q}(z)$ may be thought of as a discrete version of the Fourier transform \cite{Sinai2},
 and it is a Laurent polynomial in the variables $e^{z_1}, \dots, e^{z_d}$, 
encoding precisely all the integer points in $Q \cap \Z^d$.  In \cite{Sinai1}, we showed that for integer polytopes $\P$, the integer point transform is a complete invariant of $\P$. 
The {\bf discrete Brion formula} \cite{Brion1}
reduces the integer point transform of a rational polytope $\P$ to
a sum of the integer point transforms of the vertex tangent cones of $\P$:
\begin{equation}
\label{discrete Brion theorem}
\sigma_{\P}(z) = \sigma_{\K_{v_1}}(z)  + \cdots +  \sigma_{\K_{v_N}}(z),
\end{equation}
where the vertices of $\P$ are $v_1, \dots v_N$. For a proof of \eqref{discrete Brion theorem}, see for example \cite[Theorem 11.7]{BeckRobins}. 
It is elementary that each of the integer point transforms $\sigma_{\K_{v}}(z)$ is a rational function in the variables
$s_1:= e^{z_1}, \dots s_d:= e^{z_d}$:
\begin{equation}
\label{Elementary cone sum lemma}
\sigma_{\K_v}(z):=
\sum_{n \in \K_v \cap \Z^d}    
e^{ \langle n,  z \rangle}
= 
\frac{   \sigma_{ \Pi_v + v}(z) }
{
\prod_{k=1}^d  
\left(   1 - e^{ \langle z,  \omega_k \rangle}   \right)   
},
\end{equation} 
for each vertex $v\in \P$.
Consequently, Brion's identity  
\eqref{discrete Brion theorem}
is valid for almost all $z\in \C^d$. 
For the precise statements of the results in this paper, we formalize the meaning of a generic $z$.

\begin{definition}
\label{def:almost all z}
We call $z\in \C^d$  {\bf generic} if it belongs to the set
\begin{equation}
\mathcal G:=\{z \in \C^d \mid \ 
\langle z, w(v) \rangle 
    \not= 0, \text{ for any edge } w(v) \text{ of } \P
\}.
\end{equation}
\end{definition}
\noindent
We'll often use the words ``for almost all $z$'' to mean that $z\in \mathcal G$. We note that $\mathcal G$ is the complement of a finite, complex hyperplane arrangement.

We define  
a {\bf quasi-polynomial} to be a function 
$f(t) = \sum_{k=0}^d c_k(t) t^k$,
where each coefficient $c_k(t)$ is a periodic function of $t\in \R_{>0}$.  
It is extremely useful to define the {\bf Barnes numbers} $B_k(\va)$  to be the constant term of the Barnes polynomials:
\begin{equation}
\label{BarnesNumbers}
B_k(\va):= B_k(0, \va),
\end{equation}
valid for any $\va \in \C^d$.  This definition is the analogue of the classical $1$-variable Bernoulli numbers.

Given a half-open integer parallelepiped $\Pi$,
there is natural flat torus, defined by  $\R^d/ \Lambda$, where $\Lambda$ is the full-rank lattice generated by the vertices of $\Pi$.   We'll identify this torus with the half-open parallelepiped $\Pi$. 
We define the {\bf lattice flow} over the torus $\R^d/ \Lambda$, for any
fixed direction $v\in \R^d$, by  
\begin{equation}
\label{def:lattice flow}
   LatticeFlow(t):=  \Z^d - tv \pmod \Pi ,
\end{equation} 
as in  \S \ref{sec:polygons}, Figure \ref{Fig:lattice flow 1} below.
We will compute in a straightforward way some explicit lattice flows over  half-open parallelepipeds (see \S \ref{dimension 1} and \S \ref{sec:polygons}).  
  
\medskip
\begin{definition}
\label{discrete moments}
We define the {\bf polytope Dedekind sum}  of a half-open integer parallelepiped $\Pi \subset \R^d$ by
\begin{equation}
\label{defining sum of discrete moments}
\mu_k(\Pi, z, tv):= 
\sum_{p \in  
   \Z^d - tv \; (  \hspace{-.1cm}   \bmod \Pi) 
   }    
\langle p, z \rangle^k,
\end{equation}
for each fixed integer $k\geq 0$, 
$v \in \R^d$ and $t\in \R$.
\end{definition}
\noindent
When $t=0$, the sum \eqref{defining sum of discrete moments} is also called the {\bf $k$'th discrete moment} of $\Pi$.


\bigskip 
\begin{thm}
\label{thm:Ehrhart quasi-polynomial from Barnes}
Let $\P \subset \R^d$ be a $d$-dimensional, simple rational polytope.  
{\rm
\begin{enumerate}[(a)]
\item
\label{part a of Theorem 2}
The quasi-polynomial 
$L_\P(t):= | t\P \cap \Z^d|$ is equal to:
\begin{equation}
\label{Theorem 2, main formula}
L_\P(t) =   
  \frac{(-1)^d}{d!}    
  \sum_{v \in V}      
  \sum_{k=0}^d  
  \binom{d}{k}   
  B_k\left(t \langle v, z \rangle, \va_{v}\right)
\sum_{p \in  
    \Z^d - tv \; (  \hspace{-.1cm}   \bmod \Pi_v)
   }    
{\langle  p, z \rangle}^{d-k},
\end{equation}
valid for almost all $z \in \C^d$, and all $t>0$.

\item 
\label{part b of Theorem 2}
The coefficient of $t^r$ in $L_\P(t):= \sum_{r=0}^d c_r(t) t^r$ is equal to:
\begin{equation}
\label{quasi-coefficient, extending Linke}
c_r(t) =  \frac{(-1)^d}{d!}   
\sum_{v \in V} 
 \langle v, z \rangle^{r} 
 \sum_{j=0}^{d-r}   
   \binom{d}{r, j}
    B_{j}(\va_v)   
   \sum_{p \in  
    \Z^d - tv \; (  \hspace{-.1cm}   \bmod \Pi_v)
   }   
   {\langle  p, z \rangle}^{d-r-j}, 
\end{equation}
valid for almost all $z \in \C^d$, and all $t>0$.
The Barnes polynomials 
$B_d(x, \va)$ are defined in \eqref{Barnes}, the Barnes numbers $B_{j}(\va)$ are defined in \eqref{BarnesNumbers}, and the vector of linear forms 
 $\va_v$ is defined in \eqref{def:a_v}.
\hfill 
\hyperlink{Proof of Theorem 2}{\rm{(Proof)}}
$\square$
\end{enumerate}
}
\end{thm}

\medskip

\begin{remark}
Theorem \ref{thm:Ehrhart quasi-polynomial from Barnes}
has an extension to all rational polytopes - see Corollary \ref{cor:general rational polytopes}.
\end{remark}

\begin{remark}
It may seem rather surprising a priori that the right-hand-sides of both \eqref{Theorem 2, main formula}
and 
\eqref{quasi-coefficient, extending Linke}
appears to depend on $z \in \C^d$;  but in fact a posteriori they are independent of $z$, as the left-hand-side of 
\eqref{Theorem 2, main formula} shows.  The Barnes polynomials and the Barnes numbers are rational functions of $z$, so that locally each summand has poles in $z$;  however, their global sum over all the vertices cancels out all of these poles, due to Brion's theorem \eqref{discrete Brion theorem}.

Moreover, it is possible that judicious choices of $z$ may reduce the computational complexity of the ensuing discrete moments, and therefore of the Ehrhart polynomials and quasi-polynomials (see \S \ref{sec: further remarks}, 
Remark \ref{remark: variety idea}). 
\hfill $\hexago$
\end{remark}

\begin{remark}
We observe that for rational simple polytopes, one novelty in Theorem 
\ref{thm:Ehrhart quasi-polynomial from Barnes} is that it reduces the study of Ehrhart quasi-polynomial $L_\P(t):= | t\P \cap \Z^d |$ to the study of discrete moments of integer points in the rational parallelepipeds
$\Pi_{v}$, and to the lattice flows 
(defined in \eqref{def:lattice flow})
over its vertex parallelepipeds.  

Moreover, in the case of rational polytopes, the periodicity of the quasi-polynomail in formula 
\eqref{Theorem 2, main formula} resides only in the polytope Dedekind sums (discrete moments) 
\begin{equation}
\label{discrete moment of vertex parallelepiped}
    \mu_n(\Pi_v, z, tv):= 
\sum_{p\in\Z^d - tv \; (  \hspace{-.1cm}   \bmod \Pi_v)} 
{\langle  p, z \rangle}^{n},
\end{equation} 
for all $t>0$.  We remark that prima facie these polytope Dedekind sums may not appear to be polynomial-time computable.  However, Barvinok's algorithm \cite{Barvinok1} shows that they are (see also Question \ref{conjecture for computing discrete moments} below). 
\hfill $\hexago$
\end{remark}

\begin{remark}
We quickly recover a couple of the main results of Eva Linke's important paper \cite{Linke}.  Linke's work opened the door to a more precise study of all positive real dilates of rational polytopes.  Theorem 1.5 of  \cite{Linke} proves that the coefficients of quasi-polynomials of rational polytopes are always piecewise polynomial functions of $t$. Here it  follows easily from 
Theorem 
\ref{thm:Ehrhart quasi-polynomial from Barnes} by rewriting
\eqref{discrete moment of vertex parallelepiped} 
as follows:
\begin{align*} 
\sum_{p \in 
\Z^d - tv \; (  \hspace{-.1cm}   \bmod \Pi_v)
} 
{\langle  p, z \rangle}^{d-k}
= \sum_{q \in (\Pi_v +tv) \cap \Z^d} {\langle  q-tv, z \rangle}^{d-k} 
&=
\sum_{q \in (\Pi_v +tv) \cap \Z^d} {\Big(
\langle  q, z \rangle 
- t \langle  v, z \rangle
\Big)
}^{d-k},
\end{align*} 
which is clearly a piecewise polynomial function of $t\in \R$, because the domain of summation is piecewise constant in $t$, while the summands are growing as a polynomial of degree $d-k$ in $t$. 
Moreover, using the latter piecewise polynomial expression for the discrete moments of the parallelepiped, we may easily differentiate the $r$'th quasi-coefficient \eqref{quasi-coefficient, extending Linke}, in $t$:
\begin{align*}
&\frac{d}{dt} c_r(t) = \frac{d}{dt} \left(
\frac{(-1)^d}{d!}   \sum_{v \in V} 
 \langle v, z \rangle^{r} 
 \sum_{j=0}^{d-r}   
   \binom{d}{r, j}
    B_{j}(\va_v)   
    \sum_{q \in (\Pi_v +tv) \cap \Z^d} {\Big(
\langle  q, z \rangle 
- t \langle  v, z \rangle
\Big)
}^{d-r-j}
\right) \\
&=
-\frac{(-1)^d}{d!}   \sum_{v \in V} 
 \langle v, z \rangle^{r+1} 
 \sum_{j=0}^{d-r-1}   
   \tfrac{d!(d-r-j)}{r!j!(d-r-j)!}
    B_{j}(\va_v)   
    \sum_{q \in (\Pi_v +tv) \cap \Z^d} 
    {\Big(
\langle  q, z \rangle 
- t \langle  v, z \rangle
\Big)}^{d-r-j-1}\\
&= -(r+1) c_{r+1}(t),
\end{align*}
recovering one of the other main results in Linke's paper, 
namely  \cite[Theorem 1.6]{Linke}. As a technicality, we note that in computing the derivative above, we restrict $t$ so that 
$\Z^d - tv$ does not belong to the boundary of $\Pi_v$, allowing us to avoid any potential discontinuities.
\hfill $\hexago$
\end{remark}

\begin{remark}
Given a rational vertex $v \in \Q^d$, let's allow the lattice flow $\Z^d - tv \pmod \Pi$ to proceed for all sufficiently large positive values of $t$.   Topologically, we get a finite collection of closed geodesics on the torus, where each geodesic is an orbit of one integer point in $\Pi_v$.  As we see from Theorem \ref{thm:Ehrhart quasi-polynomial from Barnes}, the periodicity of the Ehrhart quasi-polynomials arises precisely because of the periodicity of these closed geodesics, appearing in the discrete moments of the vertex parallelepipeds of the vertices.
    \hfill $\hexago$
\end{remark}

\begin{remark}
The discrete moments \eqref{discrete moment of vertex parallelepiped} give us a very precise measure with which  the arithmetic and spatial statistics  of the distribution of integer points in each vertex parallelepiped $\Pi_v$ account for most of the complexity of Ehrhart quasi-polynomials.
   \hfill $\hexago$
\end{remark}

An immediate corollary of Theorem \ref{thm:Ehrhart quasi-polynomial from Barnes}
is the following expansion of the Ehrhart polynomial of a simple integer polytope $\P$, in terms of the Barnes polynomials and of the discrete moments of its vertex parallelepipeds. 
Its proof is trivial:  when all the vertices $v \in \P$ are integer points, and $t \in \Z$, it follows that   
$\Z^d - tv  
= \Z^d$. It follows from 
\eqref{Theorem 2, main formula}
that there is no longer any periodic behavior in $L_\P(t)$.

\begin{cor}
\label{thm:Ehrhart.from.Bernoulli.Barnes}
Let $\P \subset \R^d$ be a simple $d$-dimensional integer polytope.  Then we have:
{\rm
\begin{enumerate}[(a)]
    \item  
    \label{part a of Theorem 1}
The  Ehrhart polynomial of $\P$ is
\begin{equation}
\label{MainTheoremEhrhartPoly}
L_\P(t)  =     
\frac{(-1)^d}{d!} 
\sum_{k=0}^d  
\binom{d}{k}
\sum_{v \in V}                             
B_k(t  \langle v, z \rangle,   \va_{v})
   \sum_{p \in \Pi_v \cap \Z^d}   
   {\langle  p, z \rangle}^{d-k},
\end{equation}
valid for almost all $z \in \C^d$, 
and all $t \in \Z_{>0}$.   

\item 
 \label{part b of Theorem 1}
The coefficient of $t^r$ in 
$L_\P(t):= \sum_{r=0}^d c_r t^r$ is equal to:
\begin{equation}
\label{MainTheoremEhrhartPoly, second part}
    c_r =  \frac{(-1)^d}{d!}   
    \sum_{v \in V}                          \sum_{j=0}^{d-r}   
   \binom{d}{r, j}
    B_{j}(\va_v)  
    \langle v, z \rangle^{r}  
   \sum_{p \in \Pi_v \cap \Z^d}   
   {\langle  p, z \rangle}^{d-r-j},
\end{equation}
valid for $t \in \Z_{>0}$,
for almost all $z \in \C^d$, and for each $0\leq r \leq d$.  
The Barnes polynomials 
$B_d(x, \va)$ are defined in \eqref{Barnes}, the Barnes numbers $B_{j}(\va)$ are defined in \eqref{BarnesNumbers}, and the vector of linear forms
 $\va_v$ is defined in \eqref{def:a_v}.
\hfill 
$\square$
\end{enumerate}
}
\end{cor}

\begin{remark} 
Corollary \ref{thm:Ehrhart.from.Bernoulli.Barnes} has a direct extension to any rational polytope - 
see observation
\ref{observation:extension-principle-general-polytopes}.
\end{remark}

\begin{remark}
  We observe that Theorem \ref{thm:Ehrhart quasi-polynomial from Barnes} also implies that for an integer polytope $\P$, allowing all positive dilates $t>0$ quickly implies that the periodic components of the quasi-polynomial $L_\P(t)$ always have a period of $1$. The reason is that for integral vertices $v \in \P$, the domain of summation in the discrete moments of formula \eqref{Theorem 2, main formula}
 satisfies $\Z^d - v = \Z^d$.  It becomes an interesting question of whether there exists a smaller {\bf rational period}, for $0<t<1$,  for some integer polytopes.
\end{remark}
\bigskip \bigskip \bigskip

\begin{restatable}{thm}{Discrete moments of rational polytopes}
\label{thm:Discrete moments of rational polytopes}
Let $\P \subset \R^d$ be a $d$-dimensional, simple rational polytope. 
\begin{enumerate}[(a)]
\item 
\label{part a of discrete moments}
For each positive integer $m$, we have the following formula for the $m$'th discrete moment of positive dilates of $\P$:
\begin{align}
\sum_{p \in t\P \cap\Z^d} 
\langle p, z \rangle^m 
&= 
(-1)^d m! 
\sum_{v \in V}  
\sum_{k=0}^{d+m} 
\tfrac{1}{k!(d+m-k)!}
B_k(t  \langle v, z \rangle,   \va_{v})
\sum_{q \in 
\Z^d - tv \; (  \hspace{-.1cm}   \bmod \Pi_v)
}  
{\langle  q, z \rangle}^{d+m-k},
\end{align} 
a quasi-polynomial, valid for almost all $z\in \C^d$, and all $t>0$.
\item
The coefficient of $t^r$ in the quasi-polynomial 
$\displaystyle\sum_{p \in t\P \cap\Z^d}  \langle p, z \rangle^m:= 
\displaystyle\sum_{r=0}^{d+m} d_r(t) t^r$
of part 
\eqref{part a of discrete moments} equals:
\begin{equation}
   d_r(t) = \tfrac{(-1)^d m!}{(d+m)!}
\sum_{v \in V}  
\sum_{j=0}^{d+m-r}   
   \binom{d+m}{r, j}
    B_{j}(\va_v)  
    \langle v, z \rangle^{r}  
\sum_{q \in  \Z^d - tv \; (  \hspace{-.1cm}   \bmod \Pi_v)}   
   {\langle  q, z \rangle}^{d+m-r-j}, 
\end{equation} 
valid for almost all $z \in \C^d$, and all $t>0$.
\item 
\label{Extension of Linke's ODE}
For $0\leq r\leq d+m-1$, and for every $t>0$ away from the locally finite set of breakpoints at which at least one of the finite sets
$\bigl(\Pi_v+tv\bigr)\cap\Z^d$ changes, the quasi-polynomial coefficients $d_r(t)$ above satisfy the ordinary differential equation
\begin{align}
\label{the ODE of Linke, extended}
\frac{d}{dt}d_r(t)=-(r+1)d_{r+1}(t).
\end{align}
\end{enumerate}
The Barnes polynomials  
$B_d(x, \va)$ are defined in \eqref{Barnes},  the Barnes numbers $B_{j}(\va)$ are defined in \eqref{BarnesNumbers}, and   $\va_v$ is defined in \eqref{def:a_v}.
\hfill 
\hyperlink{Discrete moments of polytopes}{\rm{(Proof)}}
$\square$
\end{restatable}

\begin{remark}
Theorem \ref{thm:Discrete moments of rational polytopes} has a direct extension to any rational polytope - 
see observation
\ref{observation:extension-principle-general-polytopes}.
\end{remark}

\begin{remark}
The ODE \eqref{the ODE of Linke, extended} in Theorem \ref{thm:Discrete moments of rational polytopes}, part \eqref{Extension of Linke's ODE}, directly extends Linke's ODE in her seminal paper
\cite{Linke}.  The identity is understood on each open interval between successive breakpoints, or equivalently for every $t>0$ away from the locally finite set of breakpoints specified in the theorem.  It is amusing---and a bit surprising---that we obtain in
\eqref{the ODE of Linke, extended}
exactly the same ODE as in Linke's work, even though for $m>0$ our ODE describes positive {\bf discrete moments} of dilated rational polytopes. 
 \hfill $\hexago$
\end{remark}


\subsection{Vanishing identities for rational, simple polytopes}
The following collection of vanishing identities, arising from the proof of Theorem \ref{thm:Ehrhart quasi-polynomial from Barnes}, give us
 information about the {\it location of points} in simple integer polytopes.  
Hence the following identities might be helpful in integer linear programming.  From a computational perspective, these vanishing identities allow us to relate discrete moments of parallelepipeds recursively to each other, using lower discrete moments. 

\begin{thm}
\label{thm:MainVanishingIdentities}
Let $\P \subset \R^d$ be a simple $d$-dimensional rational polytope.
Then  $\P$ possesses the following set of finite vanishing identities, valid for almost all  $z\in \C^d$ and for all $t >0$:
\begin{enumerate}[(a)]
\item 
\label{part a of vanishing identities}
    For each $0\leq m \leq d-1$, we have:
\begin{align}\label{vanishing identities1}
0 
&=   
\sum_{k=0}^m 
\binom{m}{k}   
\sum_{v \in V}      
  B_k \Big(t  \langle v, z \rangle,   \va_{v} \Big)   
\sum_{q \in  \Z^d - tv \; (  \hspace{-.1cm}   \bmod \Pi_v)}       
{\langle  q, z \rangle}^{m-k}.
\end{align}  
\item 
\label{part b of vanishing identities}
Equivalently, we may express the latter identities using the Barnes
numbers, thereby removing the explicit Barnes-polynomial dependence on
$t$.  For each $(m,r)$ with $0\leq m\leq d-1$ and $0\leq r\leq m$,
we have
\begin{align}
0
&=
\sum_{v\in V}
\langle v,z\rangle^r
\sum_{j=0}^{m-r}
\frac{B_j(\va_v)}{j!(m-r-j)!}
\sum_{q\in\Z^d-tv\;(\bmod\Pi_v)}
\langle q,z\rangle^{m-r-j}.
\end{align}
The Barnes polynomials  
$B_d(x, \va)$ are defined in \eqref{Barnes},  the Barnes numbers $B_{j}(\va)$ are defined in \eqref{BarnesNumbers}, and   $\va_v$ is defined in \eqref{def:a_v}.

\hfill 
\hyperlink{Proof of Theorem 3}{\rm{(Proof)}}
$\square$
\end{enumerate}
\end{thm}

\begin{remark} 
Theorem  \ref{thm:MainVanishingIdentities} has a direct extension to any rational polytope - 
see observation
\ref{observation:extension-principle-general-polytopes}.
\end{remark}

\begin{example} 
In the special case of integer polytopes, 
 Theorem \ref{thm:MainVanishingIdentities}, part \eqref{part b of vanishing identities} with $m=1, r=0$ tells us that
\begin{align}
 \sum_{v \in \rm{V}}       \frac{ 1 }{ \prod_{k=1}^d  \langle w_k(v), z \rangle}    
\sum_{p \in \Pi_v \cap \Z^d}    \langle  p, z \rangle        
 = \frac{1}{2} \sum_{v \in \rm{V}}     \vol \Pi_{v}              
\frac{   \sum_{k=1}^d \langle w_k(v), z \rangle }{\prod_{k=1}^d \langle w_k(v), z \rangle},
\end{align}
as we show below in \eqref{independent interest, m=1, r=0}.  The latter identity sometimes allows us 
 to rapidly compute certain weighted linear combinations of discrete moments of vertex parallelepipeds.
The cases with $r =m $ appeared in \cite{BrionVergne}, while all the other cases of Theorem \ref{thm:MainVanishingIdentities} appear to be new.  We expand on these special cases $r=m$ in 
\eqref{Brion vanishing identities} below (See also \cite[Theorem 5.11]{Sinai2} for a Fourier-based proof).
\end{example}


\subsection{The constant term gives generalized reciprocity laws for polytope Dedekind sums}
The following  identities arise from the standard topological fact that the constant term of the Ehrhart polynomial of an integer polytope is always $1$.  
In all of the ensuing corollaries, we keep using the same notation for  
$B_d(x, \va)$ \eqref{BarnesNumbers},  
$B_{j}(\va)$ \eqref{Barnes}, 
and $\va_v$ \eqref{def:a_v}.

\begin{cor}[Reciprocity Law for polytope Dedekind sums]
\label{cor:ConstantCoeff}
Let $\P \subset \R^d$ be a $d$-dimensional integer, simple polytope.  Then we have:
\begin{equation} 
\label{equation weird1}    
\frac{(-1)^d}{d!} 
\sum_{k=0}^d  
\binom{d}{k}
\sum_{v \in V}                             
B_k(\va_{v})
   \sum_{p \in \Pi_v \cap \Z^d}   
   {\langle  p, z \rangle}^{d-k}=1,
\end{equation}
valid for almost all $z\in \C^d$.   
\hfill $\square$
\end{cor}
\begin{proof}
To get the constant term of the Ehrhart polynomial, we simply substitute $t=0$ in Corollary \ref{thm:Ehrhart.from.Bernoulli.Barnes} 
to obtain:
\begin{equation}
1 = L_\P(0)  =     
\frac{(-1)^d}{d!} 
\sum_{k=0}^d  
\binom{d}{k}
\sum_{v \in V}                             
B_k(0,   \va_{v})
   \sum_{p \in \Pi_v \cap \Z^d}   
   {\langle  p, z \rangle}^{d-k},
\end{equation}
and the result follows because 
$B_k(0,   \va_{v}) := B_k(\va_{v})$, the Barnes numbers. 
\end{proof}

\begin{remark}
Corollary \ref{cor:ConstantCoeff}
has a direct extension to any rational polytope - 
see observation
\ref{observation:extension-principle-general-polytopes}.
\end{remark}


\subsection{Polytope Dedekind sums, via the fractional part function}

Here we develop further the polytope Dedekind sums. 
 In particular, we now 
express each discrete moment (i.e. polytope Dedekind sum)  of the vertex parallelepiped at $v$ in terms of the fractional part function 
$\{x\}:= x - \lfloor x \rfloor$.

As before, we fix a rational, simple vertex $v$ in the rational $d$-polytope $\P$, and  we suppose that we are given its incident primitive edge vectors $w_1, \dots, w_d$, as defined in \eqref{def:vertex parallelepiped}.

\begin{lemma}[Finite abelian group formula]
\label{lemma:discrete moments via fractional parts}
For each $n\in \Z^d \cap \Pi_v$, and each $t \in \R$, we define 
\begin{equation}
\label{Lem:initial alpha_j}
\begin{pmatrix}
 \alpha_1(n, t)  \\
\vdots \\
\alpha_d(n, t)
\end{pmatrix}
=
W_v^{-1}(n-tv).
\end{equation}
Then for each fixed integer $k\geq 0$, the polytope Dedekind sum may be expressed as follows:  
\begin{equation}
\label{Lem:desired sawtooth conclusion}
\sum_{p \in  
    \Z^d - tv \; (  \hspace{-.1cm}   \bmod \Pi_v)
   }    
{\langle  p, z \rangle}^k 
= \sum_{n \in \Z^d \cap \Pi_v} 
\Big(
\big\{
\alpha_1(n, t)
\big\}
\langle w_1(v), z \rangle 
+ \cdots 
+ 
\big\{
\alpha_d(n, t)
\big\}
\langle w_d(v), z \rangle 
\Big)^k,
\end{equation}
for all $z \in \C^d$.  
\end{lemma}
\begin{proof}
By assumption, $w_1, \dots w_d$ form a basis of $\R^d$.  In particular, for any $p\in \Z^d - tv$, there are unique real numbers $\alpha_1, \dots, \alpha_d$ such that 
\begin{equation}
\label{beginning with a general point}
   p = \sum_{j=1}^d \alpha_j w_j.  
\end{equation}
 Since $p := n -tv$ for some $n \in \Z^d$, the latter sum is equivalent to 
\eqref{Lem:initial alpha_j}.
Therefore for each $n\in \Z^d \cap \Pi_v$, with its corresponding $p:= n -tv$, \eqref{beginning with a general point}
may be rewritten as 
\[
n-tv = p 
=
\sum_{j=1}^d \lfloor \alpha_j \rfloor w_j
+ 
\sum_{j=1}^d \{ \alpha_j \} w_j,
\] 
for unique $\alpha_j$'s that depend on $n, t$, and $v$.   Since any integral linear combination of the $w_j$'s by definition vanishes modulo $\Pi_v$, we reduce the latter equality mod $\Pi_v$ to obtain 
\begin{equation}
n-tv = \sum_{j=1}^d \{ \alpha_j \} w_j   
\quad ( \hspace{-.03cm} \bmod \ \Pi_v).
\end{equation}
Substituting 
$p= \displaystyle\sum_{j=1}^d \{ \alpha_j \} w_j   
\ ( \bmod \ \Pi_v)$ into the discrete moment 
$\displaystyle\sum_{p \in  
    \Z^d - tv \; (  \hspace{-.03 cm}   \bmod \Pi_v)
   }    
{\langle  p, z \rangle}^k$, now valid for all $z \in \C^d$,  gives us the desired reformulation
\eqref{Lem:desired sawtooth conclusion}. 
\end{proof}
We note that the sum on the right-hand side of  
\eqref{Lem:desired sawtooth conclusion}, considered modulo $\Pi_v$, 
may now be interpreted as a sum over a finite abelian group, namely the right-hand side of 
\eqref{Lem:desired sawtooth conclusion}.  Because one particular idea in the latter proof will be useful in some of the ensuing proofs, we record  it here, as follows.
\begin{observation}[Another description of the lattice flow]
\label{Observation: explicit lattice flow description}
Let
$\Pi_v:=W_v[0,1)^d$ be a half-open integer parallelepiped. Then the lattice flow has the following description:
\begin{equation}
\Z^d-tv \pmod{\Pi_v}
=
\left\{
W_v
\left\{
W_v^{-1}(q-tv)
\right\}
\ \mid \ 
q\in \Pi_v \cap \Z^d
\right\},
\label{eq:lattice-flow-explicit-representation}
\end{equation}
where $\{p\}:= \left(\{p_1\}, \dots, \{p_d\}\right)$.
\end{observation}
\begin{proof}
Every $q\in\Z^d$ can be written as
$q=q_0+W_vm$, with
$q_0\in\Pi_v\cap\Z^d$, and
$m\in\Z^d$.  
Therefore
\begin{equation*}
W_v\left\{W_v^{-1}(q-tv)\right\}
=
W_v\left\{W_v^{-1}(q_0-tv)+m\right\}
=
W_v\left\{W_v^{-1}(q_0-tv)\right\},
\end{equation*}
and consequently 
$
\Z^d-tv\pmod{\Pi_v}
=
\left\{
W_v\left\{W_v^{-1}(q-tv)\right\}
\ \middle|\
q\in\Pi_v\cap\Z^d
\right\}$. 
\end{proof}


\subsection{A novel formula for the codim 1 \allowbreak quasi-coefficient}

Given a $d$-dimensional rational simple polytope $\P \subset \R^d$, we recall that 
 Theorem \ref{thm:Ehrhart quasi-polynomial from Barnes} gives us explicit formulas for the quasi-polynomial
$L_\P(t) := \sum_{r=0}^d c_r(t) t^r$.
Here we give the following application of Theorem \ref{thm:Ehrhart quasi-polynomial from Barnes}, together with 
Lemma \ref{lemma:discrete moments via fractional parts} and Lemma 
\ref{lem:closed form, first discrete moment}. 

These results allow us to find a new vertex formula for  $c_{d-1}(t)$, for all $t>0$, which is  polynomial-time computable, and extends the previously known formulations for this coefficient. 

\begin{cor}
\label{Cor:novel formula for codim 1 coeff}
For each vertex $v\in V$, let $W_v$ be the edge matrix defined in
\eqref{def:edge matrix}, and define
\[
q_{v,k}
:=
\min\left\{
q\geq 1:
q e_k^T W_v^{-1}\in\Z^{1\times d}
\right\},
\qquad
\beta_{v,k}:=e_k^T W_v^{-1}v.
\]
If $\P\subset\R^d$ is a $d$-dimensional rational simple polytope,
then for every $t>0$,
\begin{equation}
\label{corrected codim one quasi coefficient}
\begin{aligned}
c_{d-1}(t)
&=
\frac{(-1)^{d-1}}{(d-1)!}
\sum_{v\in V}
\sum_{k=1}^d
\frac{
\langle v,z\rangle^{d-1}\vol(\Pi_v)
}{
\displaystyle\prod_{\substack{1\leq j\leq d\\ j\neq k}}
\langle w_j(v),z\rangle
}
\frac{1}{q_{v,k}}
\left(
\frac12-
\left\{-tq_{v,k}\beta_{v,k}\right\}
\right),
\end{aligned}
\end{equation}
for almost all $z\in\C^d$.
\hfill
\hyperlink{Codim 1 formula}{\rm{(Proof)}}
$\square$
\end{cor}

\medskip
\begin{remark} 
Corollary \ref{Cor:novel formula for codim 1 coeff}
has a direct extension to any rational polytope - 
see observation
\ref{observation:extension-principle-general-polytopes}.
\end{remark}
\begin{remark}[The Brion--Lawrence formula in codimension one]
\label{rem:Brion-Lawrence-codimension-one}

We note that if $\P\subset\R^d$ is a $d$-dimensional integer simple polytope, and we let $t$ run through the positive integers, then \eqref{corrected codim one quasi coefficient} degenerates to:
\begin{equation}
\label{corrected integral codim one coefficient}
\begin{aligned}
c_{d-1}
&=
\frac{(-1)^{d-1}}{2(d-1)!}
\sum_{v\in V}
\sum_{k=1}^d
\frac{
\langle v,z\rangle^{d-1}\vol(\Pi_v)
}{
\displaystyle
\prod_{\substack{1\leq j\leq d\\j\neq k}}
\langle w_j(v),z\rangle
}
\frac{1}{q_{v,k}}.
\end{aligned}
\end{equation}
We now explain why this special case of Corollary \ref{Cor:novel formula for codim 1 coeff}
is precisely the Brion--Lawrence volume formula
\cite[Corollary~2]{Brion1}, \cite{Lawrence}, applied separately to each of  the
facets of \(\P\).
Consider a facet \(F\) of \(\P\), and let $v$ be a vertex of $F$.   Since \(\P\) is
simple, exactly one of the edge vectors at \(v\) is not tangent to \(F\).
We denote its index by \(k=k(v,F)\).  Then the $d-1$ vectors
$\{ w_j(v) \mid   j \neq k \}$
are precisely the edge vectors of \(F\) emanating from \(v\).  We'll denote the 
lattice in the linear space parallel to 
\(F\) by 
$\Lambda_F :=
\operatorname{lin}(F-F)\cap\Z^d$.
The vector
\begin{equation}
\label{eq:primitive-normal-vector-at-facet}
u_{v,k}
:=
q_{v,k}W_v^{-T}e_k
\end{equation}
has a simple geometric meaning: it is the primitive integer normal
vector to \(F\). 
Indeed, recalling that \(W_v\) is the matrix whose \(j\)-th column is
\(w_j(v)\), and that \(q_{v,k}\) is the least positive integer for which
$
q_{v,k}W_v^{-T}e_k\in\mathbb{Z}^d
$, we have 
$w_j(v)=W_ve_j$, and therefore: 
\begin{align}
\left\langle u_{v,k},w_j(v)\right\rangle
&=
\left\langle
q_{v,k}W_v^{-T}e_k,W_ve_j
\right\rangle
\nonumber\\
&=
q_{v,k}e_k^TW_v^{-1}W_ve_j
\nonumber\\
&=
q_{v,k}e_k^Te_j
\nonumber\\
&=
q_{v,k}\delta_{kj}.
\end{align}
It follows that \(u_{v,k}\) is perpendicular to every tangent edge
\(w_j(v)\), with \(j\neq k\), while the remaining edge \(w_k(v)\) has
lattice height \(q_{v,k}\) above the affine span of \(F\).

The tangent edge vectors \(w_j(v)\), \(j\neq k\), generate a
sublattice of \(\Lambda_F\). The determinant of \(W_v\) splits into the
lattice height \(q_{v,k}\) in the normal direction and the lattice index
of the tangent edge vectors inside \(\Lambda_F\). Since
$
\vol(\Pi_v)
=
\left|\det W_v\right|
$,
we obtain
$
[
\Lambda_F:
\sum_{\substack{1\leq j\leq d\\j\neq k}}
\Z w_j(v)
]
=
\frac{\vol(\Pi_v)}{q_{v,k}}$.

So we see that  \(\vol(\Pi_v)/q_{v,k}\) is exactly the relative volume, with
respect to \(\Lambda_F\), of the vertex parallelepiped of the tangent
cone of \(F\) at \(v\).
We may therefore apply the Brion--Lawrence formula to the
\((d-1)\)-dimensional lattice polytope \(F\), using the lattice
\(\Lambda_F\). This gives
\begin{equation}
\label{eq:Brion-Lawrence-applied-to-facet}
\operatorname{vol}_{\Lambda_F}(F)
=
\frac{(-1)^{d-1}}{(d-1)!}
\sum_{v\in V(F)}
\frac{
\langle v,z\rangle^{d-1}\vol(\Pi_v)
}{
\displaystyle
\prod_{\substack{1\leq j\leq d\\
j\neq k(v,F)}}
\langle w_j(v),z\rangle
}
\frac{1}{q_{v,k(v,F)}}.
\end{equation}
Finally, we have:
\begin{align}
c_{d-1}
&=
\label{classical codim 1 formula}
\frac12
\sum_{\substack{F\subset\P\\ \dim F=d-1}}
\operatorname{vol}_{\Lambda_F}(F) \\
&=
\frac12
\sum_{\substack{F\subset\P\\ \dim F=d-1}}
\frac{(-1)^{d-1}}{(d-1)!}
\sum_{v\in V(F)}
\frac{
\langle v,z\rangle^{d-1}\vol(\Pi_v)
}{
\displaystyle
\prod_{\substack{1\leq j\leq d\\
j\neq k(v,F)}}
\langle w_j(v),z\rangle
}
\frac{1}{q_{v,k(v,F)}} \\
&=
\frac{(-1)^{d-1}}{2(d-1)!}
\sum_{v\in V}
\sum_{k=1}^d
\frac{
\langle v,z\rangle^{d-1}\vol(\Pi_v)
}{
\displaystyle
\prod_{\substack{1\leq j\leq d\\j\neq k}}
\langle w_j(v),z\rangle
}
\frac{1}{q_{v,k}}.
\label{last line of codim 1 recap}
\end{align}
where we used the classical formula for the codimension-one Ehrhart
coefficient in 
\eqref{classical codim 1 formula}, and we 
reindexed the resulting vertex--facet incidences by the pairs
\((v,k)\) in \eqref{last line of codim 1 recap}. 
 \hfill $\hexago$
\end{remark}


\subsection{Open polytopes, and a reciprocity law for discrete moments}
\label{sec:moment-reciprocity}

Here we let 
$\mathcal{P}^{o}$ be the interior of $\mathcal P$.  In this section, we extend Ehrhart's well-known reciprocity from lattice-point enumeration to the
discrete moments considered in Theorem~2. 
  For
$m\in\mathbb{Z}_{\geq 0}$, $t>0$, and $z\in\mathbb{C}^d$, define the
closed and open moments by
\begin{equation}
\begin{aligned}
M_m(\mathcal{P};t,z)
&:=
\sum_{p\in t\mathcal{P}\cap\mathbb{Z}^d}
\langle p,z\rangle^m,\\
M_m^{o}(\mathcal{P};t,z)
&:=
\sum_{p\in t\mathcal{P}^{o}\cap\mathbb{Z}^d}
\langle p,z\rangle^m.
\end{aligned}
\label{eq:closed-open-discrete-moments}
\end{equation}

For any $s\in\mathbb{R}$ and generic $z\in\mathbb{C}^d$, we continue to use the same notation for the extension of our explicit $m$th moment by
\begin{equation}
\begin{aligned}
M_m(\mathcal{P};s,z)
:=
(-1)^d m!
\sum_{v\in V}
\sum_{k=0}^{d+m}
\frac{
B_k\!\left(s\langle v,z\rangle,\mathbf{a}_v\right)
}{
k!\,(d+m-k)!
}
\sum_{q\in(\mathbb{Z}^d-sv)\cap\Pi_v}
\langle q,z\rangle^{d+m-k}.
\end{aligned}
\label{eq:barnes-extension-of-moments}
\end{equation}
Equivalently, using the lattice-flow notation of Theorem~2, the final
sum in \eqref{eq:barnes-extension-of-moments} is taken over
\[
q\in\mathbb{Z}^d-sv\,(\operatorname{mod}\Pi_v).
\]
We note that for negative $s$, 
$M_m(\mathcal{P};s,z)$ refers to the evaluation of 
formula \eqref{eq:barnes-extension-of-moments}, and not to a moment
sum over the reflected polytope $s\mathcal{P}$.

\begin{thm}[Ehrhart-type Reciprocity for discrete moments]
\label{thm:moment-reciprocity}
Let $\mathcal{P}\subset\R^d$ be a $d$-dimensional 
rational polytope.  Then, for every
$m\in\Z_{\geq 0}$ and every $t>0$:
\begin{equation}
M_m^{o}(\mathcal{P};t,z)
=
(-1)^{d+m}
M_m(\mathcal{P};-t,z),
\label{eq:moment-reciprocity-statement}
\end{equation}
for generic $z\in\mathbb{C}^d$.  
The identity extends to every
$z\in\mathbb{C}^d$ by polynomial continuation.

\noindent
In particular, we note that the case $m=0$ gives us the classical Ehrhart reciprocity law, for all postive dilates:
\begin{equation}
\label{Ehrhart reciprocity}
L_{\mathcal{P}^{o}}(t)
=
(-1)^d L_{\mathcal{P}}(-t).
\end{equation}
\hfill 
\hyperlink{Proof of Theorem: moment-reciprocity}
{\rm{(Proof)}}
$\square$
\end{thm}

\noindent
An equivalent statement of Theorem \ref{thm:moment-reciprocity} is:
\begin{equation}
\begin{aligned}
\sum_{p\in t\mathcal{P}^{o}\cap\mathbb{Z}^d}
\langle p,z\rangle^m
&=
(-1)^{d+m}
M_m(\mathcal{P};-t,z).
\end{aligned}
\label{eq:moment-reciprocity-expanded-statement}
\end{equation}


From Theorem \ref{thm:Ehrhart quasi-polynomial from Barnes}
and Ehrhart reciprocity 
\eqref{Ehrhart reciprocity}, we immediately obtain the following explicit formula for open polytopes. 

\begin{cor}
\label{Ehrhart for open polytopes} 
Let $\mathcal{P}\subset\mathbb{R}^d$ be a $d$-dimensional  simple rational polytope, and
let $\mathcal{P}^{o}$ denote its interior. Then the following statements hold.

\begin{enumerate}[(a)]

\item 
The quasi-polynomial of $\mathcal{P}^{o}$ is equal to
\begin{equation}
L_{\mathcal{P}^{o}}(t)
=
\frac{1}{d!}
\sum_{v\in V}
\sum_{k=0}^{d}
\binom{d}{k}
B_k\!\left(-t\langle v,z\rangle,\mathbf{a}_v\right)
\sum_{p\in\mathbb{Z}^d+tv\,(\operatorname{mod}\Pi_v)}
\langle p,z\rangle^{d-k},
\label{eq:cor4-open-ehrhart-quasipolynomial}
\end{equation}
valid for almost all $z\in\mathbb{C}^d$, and all $t>0$.

\item
Writing
\[
L_{\mathcal{P}^{o}}(t)
=
\sum_{\ell=0}^{d}b_{\ell}(t)t^{\ell},
\]
the coefficient $b_r(t)$ of $t^r$, for $0\leq r\leq d$, is equal to
\begin{equation}
b_r(t)
=
\frac{(-1)^r}{d!}
\sum_{v\in V}
\langle v,z\rangle^r
\sum_{j=0}^{d-r}
\binom{d}{r,j}
B_j(\mathbf{a}_v)
\sum_{p\in\mathbb{Z}^d+tv\,(\operatorname{mod}\Pi_v)}
\langle p,z\rangle^{d-r-j},
\label{eq:cor4-open-ehrhart-quasi-coefficients}
\end{equation}
valid for almost all $z\in\mathbb{C}^d$, and all $t>0$.
\end{enumerate}
\end{cor}

\begin{remark} 
Corollary 
\ref{Ehrhart for open polytopes} 
has a direct extension to any rational polytope - 
see observation
\ref{observation:extension-principle-general-polytopes}.
\end{remark}

\begin{remark}
We notice that the only difference between the periodicity of the quasi-coefficients for the interior of $\P$, versus the corresponding periodicity of the quasi-coefficients for $\P$, is that for each vertex the lattice flow on the torus (defined in \eqref{def:lattice flow} above)  is now flowing in the opposite direction.  
 \hfill $\hexago$
\end{remark}

A lattice polytope is called {\bf hollow} if its interior contains no integer points. There is considerable current research on classifying  hollow simplices \cite{Averkov.etal}. 
An immediate consequence of the Corollary \ref{Ehrhart for open polytopes} is the following equivalence for this property.

\begin{cor}
\label{cor:hollow polytopes}
Let $\P \subset \R^d$ be a $d$-dimensional integer simple polytope.  Then $\P$ is hollow if and only if 
\begin{equation}
0=  \sum_{v \in V}      
  \sum_{k=0}^d  
  \binom{d}{k}   
  B_k\big(
  -  \langle v, z \rangle,   \va_{v}
  \big)
  \sum_{p \in \Pi_v \cap Z^d}   
   {\langle  p, z \rangle}^{d-k},
\end{equation}
for almost all $z \in \C^d$.
\hfill $\square$
\end{cor}

\begin{remark} 
Corollary 
\ref{cor:hollow polytopes} 
has a direct extension to any rational polytope - 
see observation
\ref{observation:extension-principle-general-polytopes}.
\end{remark}


\subsection{General rational polytopes}
Here we prove that the formulas of Theorem 
\ref{thm:Ehrhart quasi-polynomial from Barnes} extend to any rational polytope, by employing a theoretically easy triangulation argument.   Similarly, all of the results in this paper extend to all rational polytopes, as we show in observation \ref{observation:extension-principle-general-polytopes} below, except of course for results on smooth polytopes.
 
Suppose we are given any $d$-dimensional rational polytope $\P \subset \R^d$.  Given a vertex tangent cone $\K_v$ of $\P$,  we can triangulate $\K_v$ 
into $N_v$ partially-open, disjoint, rational simplicial subcones:
\begin{equation*}
\K_v = \bigcup_{j=1}^{N_v} \K_v(j).
\end{equation*}
 In such a triangulation, none of the latter distinct simplicial subcones overlap at all, and this comes at a cost of using partially-open simplicial subcones.   The only other additional idea that we need here is the following standard fact about the vertex parallelepiped of a partially-open simplicial cone. 

\begin{lemma}
\label{lem:extension to any rational polytope}
Suppose we are given a partially-open
$d$-dimensional simplicial, integer cone $\K\subset \R^d$, with primitive edge vectors $w_1, \dots w_d$, defined by:
\begin{equation}
    \K:= \left\{
\sum_{j=1}^d \lambda_j w_j \ \mid \
0 < \lambda_1, \dots,  \lambda_m, 
\text{ and }
0 \leq \lambda_{m+1}, \dots,  \lambda_d 
\right\}.
\end{equation}
Then:
\begin{enumerate}  [(a)]
\item 
Its vertex parallelepiped has the description: 
\begin{equation}
\label{def:vertex parallelepiped for partially open cone}
\Pi
:=
\left\{
\sum_{j=1}^d \theta_j w_j
\;\middle|\;
0 < \theta_1, \dots, \theta_m \le 1,\;
\text{ and }
0 \leq \theta_{m+1},\;
\dots, \theta_d < 1
\right\}.
\end{equation}
\item We also have:
\begin{equation*}
\sigma_{\K}(z)=
\frac{   \sigma_{ \Pi}(z) }
{
\prod_{k=1}^d  
\left(   1 - e^{ \langle z,  \omega_k \rangle}   \right)   
},
\end{equation*} 
    where  $\sigma_{ \Pi}(z)$ was defined more generally in \eqref{def: general def of integer point transform}.  
\end{enumerate}
\hfill $\square$
\end{lemma}
In other words, for each facet that is excluded from the closure of $\K$, the ``opposite'' facet of its vertex parallelepiped $\Pi$ is included.  For a proof of Lemma \ref{lem:extension to any rational polytope}, see for example \cite{Barvinok2}. 

This variant of a fundamental domain, defined for any partially-open simplicial cone $\K$, also tiles $\K$ by translations with the same edge vectors $w_j$.  For each simplicial subcone $\K_{v,j}$, let
\begin{equation}
\label{def:subcone edge matrix}
W_{v,j}
:=
\begin{bmatrix}
w_1(v,j) & \cdots & w_d(v,j)
\end{bmatrix},
\qquad
\va_{v,j}:=W_{v,j}^T z
=
\begin{pmatrix}
\langle w_1(v,j),z\rangle\\
\vdots\\
\langle w_d(v,j),z\rangle
\end{pmatrix}
\in\C^d,
\end{equation}
where $w_1(v,j),\dots,w_d(v,j)$ are the primitive edge vectors of $\K_{v,j}$. 
We let $\Pi_{v,j}$ be the vertex parallelepiped of $\K_{v,j}$, defined in \eqref{def:vertex parallelepiped for partially open cone}.  As always, the Barnes polynomials  
$B_d(x, \va)$ are defined in \eqref{Barnes},  the Barnes numbers $B_{j}(\va)$ are defined in \eqref{BarnesNumbers}, and   $\va_v$ is defined in \eqref{def:a_v}.


\begin{cor}
\label{cor:general rational polytopes}
Let $\P \subset \R^d$ be a $d$-dimensional rational polytope.  Then using the notation above, the quasi-polynomial 
$L_\P(t):= | t\P \cap Z^d|$ is  equal to:
\begin{equation}
L_\P(t) =   
  \frac{(-1)^d}{d!}    
  \sum_{v \in V}   
  \sum_{j=1}^{N_v}
  \sum_{k=0}^d  
  \binom{d}{k}   
B_k\left(t \langle v, z \rangle, \va_{v, j}\right)
\sum_{p \in  
    \Z^d - tv \; (  \hspace{-.1cm}   \bmod \Pi_{v, j})
   }    
{\langle  p, z \rangle}^{d-k},
\end{equation}
valid for almost all $z \in \C^d$, and all $t>0$.
\hfill 
\hyperlink{Proof of Corollary:extension of main theorem to rational polytopes}
{\rm{(Proof)}}
$\square$
\end{cor}

\begin{observation} [Extension principle for arbitrary rational polytopes]
\label{observation:extension-principle-general-polytopes}
Although several of the preceding results are stated for simple rational
polytopes, the simplicity assumption was only used to ensure that every
vertex tangent cone is simplicial.  Given any rational polytope $\P\subset \R^d$, the 
following argument shows that a triangulation of the vertex tangent cones of $\P$ gives us a direct extension of most of the results of this paper. 

For each vertex \(v\in \P\), 
we fix a disjoint
half-open  simplicial decomposition of its vertex tangent cone $C_v$:
\begin{equation}
\label{eq:general-extension-half-open-decomposition}
C_v
=
\bigsqcup_{j=1}^{N_v} C_{v,j},
\end{equation}
as in the proof of Corollary
\ref{cor:general rational polytopes}. We let
$
w_1(v,j),\ldots,w_d(v,j),
\ W_{v,j},
\ \Pi_{v,j},
\ \va_{v,j}$
denote the corresponding primitive ray generators, edge matrix,
half-open fundamental parallelepiped, and vector of linear forms, respectively.

Every identity in this paper whose proof uses simplicity only through
the simpliciality of the vertex tangent cones extends to arbitrary rational
polytopes after replacing each vertex contribution by the sum of the
contributions of its simplicial subcones (except of course the results for smooth polytopes). More precisely, the replacement
has the form
\begin{equation}
\label{eq:general-extension-replacement-rule}
\sum_{v\in V}
F\bigl(v,W_v,\Pi_v,\va_v\bigr)
\quad\longmapsto\quad
\sum_{v\in V}\sum_{j=1}^{N_v}
F\bigl(v,W_{v,j},\Pi_{v,j},\va_{v,j}\bigr).
\end{equation}
No additional correction terms or inclusion--exclusion signs occur,
because the decomposition in
\eqref{eq:general-extension-half-open-decomposition} is disjoint.

Consequently, our formulas for the Ehrhart quasi-polynomial and
its quasi-coefficients, the corresponding formulas for discrete moments
and their quasi-coefficients, the Laurent-coefficient vanishing identities,
and the constant-term identities all possess analogous formulations for
arbitrary rational polytopes. The same principle also applies to the formulas
for open polytopes and reciprocity, using the corresponding complementary
half-open decompositions of the open tangent cones.

The individual summands in
\eqref{eq:general-extension-replacement-rule} may depend on the chosen
triangulations and generic perturbations, but their total does not, since
it equals the integer-point transform of the original tangent cone. The
genericity condition on \(z\) is understood to exclude the finite collection
of hyperplanes
$
\left\{
z\in\C^d:
\langle w_i(v,j),z\rangle=0
\right\}$, where 
$v\in V(\P),
1\leq j\leq N_v$, 
and
$1\leq i\leq d$.
\hfill $\hexago$
\end{observation}


\subsection{Smooth polytopes}

A {\bf smooth polytope} is a simple, integer $d$-dimensional polytope $\P\subset \R^d$ such that at each vertex $v \in \P$, its
primitive edge vectors
 $w_1(v), \dots w_d(v)$ form a basis for $\Z^d$. 
The term ``smooth polytopes'' comes from toric geometry, and smooth polytopes are also called ``totally unimodular polytopes" in integer programming, and ``delzant polytopes" in symplectic geometry.  
It is well known that on the one hand, it is easy to prove numerous results regarding the discrete geometry of smooth polytopes.  From the present perspective, this ease is due to the fact that most of the discrete moments at their vertices vanish, as we observe next, leading to an extreme simplification of the main results.  On the other hand, there are many open problems for smooth polytopes. 

First, specializing Theorem \ref{thm:Ehrhart quasi-polynomial from Barnes} to smooth polytopes immediately gives us the following formula for their Ehrhart polynomials.  The proofs of the rest of the corollaries are omitted, because their immediate proofs are identical to that of  
Corollary \ref{Ehrhart for smooth polytopes}.  
\begin{cor}
\label{Ehrhart for smooth polytopes}
Let $\P \subset \R^d$ be a $d$-dimensional smooth polytope.   Then we have:
\begin{equation*}
L_\P(t)  
=
\frac{(-1)^d}{d!}    
\sum_{v \in V}    
  B_d(t  \langle v, z \rangle, \va_{v}),
\end{equation*}
for all $t \in \Z_{>0}$, and almost all $z \in \C^d$. Equivalently, we have:
\begin{equation}
L_\P(t)  
=
\frac{(-1)^d}{d!}  
\sum_{j=0}^d  
\left(
\binom{d}{j}
\sum_{v \in V}      
B_{j}(\va_v) \langle v, z \rangle^{d-j}
\right)
t^{d-j},
\end{equation}
for all $t \in \Z_{>0}$, and almost all $z \in \C^d$.
\hfill 
\hyperlink{Proof of smooth formula}{\rm{(Proof)}}
$\square$
\end{cor}
\noindent 
Specializing Theorem \ref{thm:Discrete moments of rational polytopes}
to smooth polytopes gives the following moment formulas for smooth polytopes.
\begin{cor}
\label{cor:Discrete moments of rational smooth polytopes}
Let $\P \subset \R^d$ be a $d$-dimensional smooth polytope. 
For each positive integer $m$, we have the following formula for the $m$'th discrete moment of positive dilates of $\P$:
\begin{align}
\sum_{p \in t\P \cap\Z^d} 
\langle p, z \rangle^m 
&= 
(-1)^d 
\frac{m!}{(d+m)!}
\sum_{v \in V}  
B_{d+m}(t  \langle v, z \rangle,   \va_{v})
\end{align} 
valid for almost all $z\in \C^d$, and all $t \in \Z_{>0}$.
\hfill $\square$
\end{cor}

\noindent 
Specializing Theorem \ref{thm:MainVanishingIdentities} to smooth polytopes immediately gives the following set of vanishing identities. 
\begin{cor}
Let $\P \subset \R^d$ be a $d$-dimensional smooth polytope.
Then  $\P$ possesses the following set of finite vanishing identities, valid for almost all  $z\in \C^d$, and for all $t \in \Z_{>0}$:
\begin{enumerate}[(a)]
\item 
\label{part a of vanishing identities for smooth polytopes}
    For each $0\leq m \leq d-1$, we have:
\begin{align}
0 &=     
\sum_{v \in V}      
  B_m \Big(t  \langle v, z \rangle,   \va_{v} \Big).
\end{align}  
\item 
\label{part b of vanishing smooth identities}
Equivalently, we may also express the latter set of identities by using the Barnes numbers, and removing the variable $t$.  For each $(m, r)$ with  $0\leq m \leq  d-1,   \text{  and }  0\leq r \leq m$, and for almost all $z\in \C^d$ we have:
\begin{align}
0 
&=  
\sum_{v \in \rm{V}}    
   \langle v, z \rangle^{r}       
    B_{m-r}(\va_{v}).
\end{align}  
\hfill $\square$
\end{enumerate}
\end{cor}

\noindent
Specializing Corollary \ref{cor:ConstantCoeff} to smooth polytopes, we have the following.
\begin{cor}
\label{constant term 1 identity for smooth polytopes}
Let $\P \subset \R^d$ be a $d$-dimensional smooth polytope.  Then we have:
\begin{equation} 
\label{amusing identities for constant term 1}
1=     
\frac{(-1)^d}{d!}  
\sum_{v \in V}                             
B_d(\va_{v}),
\end{equation}
for almost all  $z\in \C^d$.
\hfill $\square$
\end{cor}

For example, it follows from Corollary \ref{constant term 1 identity for smooth polytopes}
that for integer smooth polygons,  we have the amusing identity:
\begin{equation}
12 =   
 \sum_{v \in \rm{V}}  
 \left(    \frac{\langle w_2(v), z \rangle}{\langle w_1(v), z \rangle} +
\frac{\langle w_1(v), z \rangle}
{\langle w_2(v), z \rangle}   +  3
\right),
\end{equation}
valid for almost all $z \in \C^2$ (see Example  \ref{first example}).   The integer $12$ is easy to explain - it arises naturally from the second classical Bernoulli number $B_2 = \frac{1}{6}$. Similarly, Corollary \ref{constant term 1 identity for smooth polytopes} gives us  a corresponding identity for all $3$-dimensional smooth polytopes:
 \begin{align*} 
24  =   \sum_{v \in \rm{V}}     
 \left(
\tfrac{\langle w_1(v), z \rangle}{\langle w_2(v), z \rangle}+\tfrac{\langle w_1(v), z \rangle}{\langle w_3(v), z \rangle}+
\tfrac{\langle w_2(v), z \rangle}{\langle w_1(v), z \rangle}+\tfrac{\langle w_2(v), z \rangle}{\langle w_3(v), z \rangle}+\tfrac{\langle w_3(v), z \rangle}{\langle w_1(v), z \rangle}+\tfrac{\langle w_3(v), z \rangle}{\langle w_2(v), z \rangle}
  +3
 \right),
\end{align*}
valid for almost all $z \in \C^3$ (see Example \ref{Example: dimension 3, the 24-identity}).
Such identities appear to be new.


\begin{remark}[A shortened form of Theorem~\ref{thm:Ehrhart quasi-polynomial from Barnes}, part \eqref{part a of Theorem 2}]
The Barnes
polynomials satisfy the Appell translation identity
$B_d(X+Y,\bm a)
=
\sum_{k=0}^{d}
\binom{d}{k}
B_k(X,\bm a)\,Y^{d-k}$.
For a lattice point
\begin{equation*}
\bm q\in
\bigl(t\bm v+\Pi_{\bm v}\bigr)\cap\mathbb Z^d,
\label{eq:q-in-translated-parallelepiped}
\end{equation*}
we take
$X=t\langle \bm v,\bm z\rangle,
\
Y=\langle \bm q-t\bm v,\bm z\rangle$. 
So we have 
$X+Y
=
t\langle \bm v,\bm z\rangle
+
\langle \bm q-t\bm v,\bm z\rangle
=
\langle \bm q,\bm z\rangle$. 
Consequently,  
Theorem~\ref{thm:Ehrhart quasi-polynomial from Barnes}, part \eqref{part a of Theorem 2}
may equivalently
be written in the shortened form
\begin{equation}
\label{eq:shortened-form-theorem-1a}
L_P(t)
=
\frac{(-1)^d}{d!}
\sum_{\bm v\in V(P)}
\;
\sum_{\bm q\in
\left(t\bm v+\Pi_{\bm v}\right)\cap\mathbb Z^d}
B_d\bigl(\langle \bm q,\bm z\rangle,\bm a_{\bm v}\bigr).
\end{equation}
Although \eqref{eq:shortened-form-theorem-1a} is more compact, the
 form of 
 Theorem~\ref{thm:Ehrhart quasi-polynomial from Barnes}, part \eqref{part a of Theorem 2}
 is generally preferable because 
it isolates the purely periodic part as a
discrete moment over a translated fundamental parallelepiped.
Moreover, these discrete moments are natural higher-dimensional
extensions of classical Dedekind sums, and refer to them as
\emph{polytope Dedekind sums}. 
\end{remark}

The literature on Ehrhart polynomials and quasi-polynomials is vast, and the reader may consult, for example,  the books
\cite{Barvinok2}  \cite{BeckRobins} \cite{HibiBook} 
\cite{LasserreBook}
\cite{Sinai2} \cite{Stanleybook},
and  papers 
\cite{Baldoni.etal}
\cite{Baldoni.etc}
\cite{Barany}
\cite{Barvinok1} 
\cite{BeckSamWoods}
\cite{BrionVergne}
\cite{DiazRobins}
\cite{DiazQuangRobins}
\cite{Ehrhart1}
\cite{StavrosPommersheim}
\cite{Godinho.etal}
\cite{HasseNillPayne}
\cite{Henk1}
\cite{HenkLinke}
\cite{HenkSchurmannWills}
\cite{Hibi}
\cite{Jochemko.etal}
\cite{KarshonSternbergWeitsman}
\cite{LagariasZiegler}
\cite{LasserreZeron}
\cite{Macdonald1} 
\cite{Macdonald2}
\cite{McAllisterMoriarity}
\cite{MeszarosMorales}
\cite{Morelli} 
\cite{Pommersheim1} 
\cite{Reznick}
\cite{LudwigSilverstein}
\cite{Stanley2}.


\bigskip 
\section{A formula for the first discrete moment}
We continue to develop the polytope Dedekind sums (aka discrete moments) of integer parallelepipeds, by giving a useful evaluation of 
$\mu_1(\Pi, z, tv)$.    
We evaluate this first discrete moment by using Lemma  \ref{lemma:discrete moments via fractional parts}.  

\begin{lemma}
\label{lem:closed form, first discrete moment}
The first
discrete moment of  the half-open integer parallelepiped $\Pi_v$ has the following evaluation:
\begin{equation}
\label{closed form for first discrete moment}
\sum_{p\in \Z^d-tv\;(\bmod \Pi_v)}
\langle p,z\rangle
=
\vol(\Pi_v)
\sum_{k=1}^d
\frac{1}{q_k}
\left(
\frac{q_k-1}{2}
+
\left\{-tq_k\beta_k\right\}
\right)
\langle w_k(v),z\rangle,
\end{equation}
where
\[
q_k
:=
\min\left\{
q\geq 1:
q e_k^T W_v^{-1}\in\Z^{1\times d}
\right\},
\qquad
\beta_k:=e_k^T W_v^{-1}v,
\]
and $z\in\C^d$.
\end{lemma}
\begin{proof}
Recall that $\Lambda_v=W_v\Z^d$ by \eqref{def:edge lattice}, and let
\[
G_v:=\Z^d/\Lambda_v.
\]
It is a standard fact from the geometry of numbers that
$
|G_v|=[\Z^d:\Lambda_v]=\vol(\Pi_v).
$
By Lemma~\ref{lemma:discrete moments via fractional parts}, the
$k$th coordinate of the coefficient vector associated with
$n-tv$ is
\[
\alpha_k(n,t)=e_k^T W_v^{-1}(n-tv).
\]
The map  $\psi: n\longmapsto W_v^{-1}n\pmod{\Z^d}$
identifies $\Z^d\cap\Pi_v$ with the finite group
\[
W_v^{-1}\Z^d/\Z^d\cong \Z^d/\Lambda_v.
\]
The $k$th coordinate of $\psi$ defines a group homomorphism
\[
\varphi_k:G_v\longrightarrow \R/\Z,
\qquad
\varphi_k(n):=e_k^T W_v^{-1}n\pmod 1.
\]
Its image is the cyclic subgroup of $\R/\Z$ of order $q_k$.  Each
point of $\operatorname{im}(\varphi_k)$ therefore occurs exactly
$|G_v|/q_k$ times as $n$ ranges over $\Z^d\cap\Pi_v$.  Since
$\beta_k=e_k^T W_v^{-1}v$, we have
\[
\left\{\alpha_k(n,t)\right\}
=
\left\{\varphi_k(n)-t\beta_k\right\},
\]
and therefore
\[
\sum_{n\in\Z^d\cap\Pi_v}
\left\{\alpha_k(n,t)\right\}
=
\frac{|G_v|}{q_k}
\sum_{r=0}^{q_k-1}
\left\{
\frac{r}{q_k}-t\beta_k
\right\}.
\]
We use the elementary finite-cyclic-sum identity
$
\sum_{r=0}^{q-1}
\left\{
\frac rq+x
\right\}
=
\frac{q-1}{2}+\{qx\}$. 
Applying it with $x=-t\beta_k$ gives
\begin{equation}
\label{sum of alpha fractional parts}
\sum_{n\in\Z^d\cap\Pi_v}
\left\{\alpha_k(n,t)\right\}
=
\frac{|G_v|}{q_k}
\left(
\frac{q_k-1}{2}
+
\left\{-q_kt\beta_k\right\}
\right).
\end{equation}
Finally, Lemma~\ref{lemma:discrete moments via fractional parts} and
\eqref{sum of alpha fractional parts} together yield
\begin{align*}
\sum_{p\in\Z^d-tv\;(\bmod\Pi_v)}
\langle p,z\rangle
&=
\sum_{n\in\Z^d\cap\Pi_v}
\sum_{k=1}^d
\left\{\alpha_k(n,t)\right\}
\langle w_k(v),z\rangle\\
&=
\vol(\Pi_v)
\sum_{k=1}^d
\frac{1}{q_k}
\left(
\frac{q_k-1}{2}
+
\left\{-q_kt\beta_k\right\}
\right)
\langle w_k(v),z\rangle.
\end{align*}
\end{proof}


\section{Small dimensions}
\label{sec:small dimensions}

In this section we show how to specialize the formula of 
Theorem \ref{thm:Ehrhart quasi-polynomial from Barnes} 
to obtain the classical Ehrhart polynomial results in dimensions $1$ and $2$.  

\subsection{Dimension 1: intervals}
\label{dimension 1}

It is instructive to first consider the obligatory $1$-dimensional examples. 
\begin{example}  
Consider the $1$-dimensional integer polytope $\P:= [a, b]$, where $a, b \in \Z$. 
Here the vertex parallelepipeds
at both vertex $a$ and vertex $b$ are unit intervals.  
The relevant Barnes polynomials in dimensions $1$ are 
$B_0(t, a_1) = \frac{1}{a_1}$, and 
$ B_1(t, a_1) = \frac{1}{a_1} t - \frac{1}{2}$
(see \eqref{B_0(t, vector a} and \eqref{B_1(t, (a_1))}).
Here the vector of linear forms is just $\va_{v} = z$ for the vertex $v:=a$, and
$\va_{v} = -z$ for the vertex $v:=b$.
 In dimension $1$ all polytopes are smooth, so Corollary \ref{Ehrhart for smooth polytopes} gives us: 
\begin{align*}
L_\P(t) 
&=
-\sum_{v \in \rm{V}}   
B_1\Big(       
t  \langle v, z \rangle,   \va_{v} 
\Big) =
-1\cdot B_1(taz, z) - 1\cdot B_1(tbz, -z) \\
&= -\Big(   \frac{taz}{z} - \frac{1}{2}  + \frac{tbz}{-z} - \frac{1}{2}  \Big) 
= t(b-a) + 1,
\end{align*}
recovering the simplest $1$-dimensional case.
\hfill $\square$
\end{example}

\begin{example}
More generally, consider the one-dimensional interval
$\P=[\alpha,\beta]$, where $\alpha,\beta\in\R$. The vertex parallelepipeds are
$\Pi_\alpha=[0,1)$ and 
$\Pi_\beta=(-1,0]$. 
For rational $\alpha$ and $\beta$, Theorem~\ref{thm:Ehrhart quasi-polynomial from Barnes}
gives the following Barnes-polynomial derivation. 
The representatives of the two lattice flows are 
$\{-t\alpha\}\in[0,1)$ and $-\{t\beta\}\in(-1,0]$.  Therefore
\begin{align*}
\left|t\P\cap\Z\right|
&=
-\sum_{v\in V}
B_1\!\left(t\langle v,z\rangle,\va_v\right)
-
\sum_{v\in V}
B_0\!\left(t\langle v,z\rangle,\va_v\right)
\sum_{p\in\Z-tv\;(\bmod\Pi_v)}
\langle p,z\rangle\\
&=
-\left(
\frac{t\alpha z}{z}-\frac12
+
\frac{t\beta z}{-z}-\frac12
\right)
-
\left(
\frac{1}{z}\{-t\alpha\}z
+
\frac{1}{-z}\bigl(-\{t\beta\}z\bigr)
\right)\\
&=
 t(\beta-\alpha)
-\{t\beta\}
-\{-t\alpha\}
+1,
\end{align*}
recovering the known formula.
\hfill $\square$
\end{example}


\subsection{Dimension 2: integer polygons}
\label{sec:polygons}

Suppose we are given any integer polygon $\P$, defined as the convex hull of its vertices $\{ v_1, \cdots v_N \}$.  At each vertex $v$, we are given its primitive edge vectors $w_1(v), w_2(v)$.
We recall that the only Barnes polynomials that occur in dimension $2$ are:
\begin{equation}
\label{ex: dim 2, B_0, B_1}
 B_0(t, \icol{ a_1 \\ a_2}) = \frac{1}{a_1  a_2}, 
 \quad \quad
 B_1(t, \icol{ a_1 \\ a_2}) 
 = \frac{1}{a_1  a_2} t -\frac{1}{2}   
 \frac{a_1 + a_2}{a_1 a_2},   
 \end{equation}
and
\begin{align}
\label{ex: dim 2, B_2}
B_2\Big(t, \icol{ a_1 \\ a_2}\Big) &= 
t^2 \left(   \frac{1}{a_1 a_2}  \right)
- t \left(   \frac{1}{a_2}   +  \frac{1}{a_1}            \right) 
+ \frac{1}{6}\left(    \frac{a_2}{a_1} + \frac{a_1}{a_2}      
\right) + \frac{1}{2},
\end{align}
using \eqref{B_2(t, (a1, a2))} of 
Appendix \ref{Basic properties of Barnes}.
We have the following
edge vectors for the vertex tangent cones at each of the vertices $v_1, v_2, v_3$:
\[ 
v_1:= \icol{0\\0}  \implies 
\omega_1(v_1)= \icol{1\\0}, \omega_2(v_1) = \icol{0\\1}.
\]
\[
v_2:= \icol{1\\0}  \implies 
\omega_1(v_2)= \icol{-1   \\     \  1}, \omega_2(v_2) = \icol{ -1 \\ \   0}.  
\]
\[
v_3:= \icol{0\\1} \implies   
\omega_1(v_3)= \icol{  \   1  \\  - 1}, \omega_2(v_3) = \icol{ \  0 \\ -1}.
\]
We claim that an explicit description for the Ehrhart polynomial of  
$\P$, using Corollary \ref{thm:Ehrhart.from.Bernoulli.Barnes}, is as follows:
\begin{align}
\label{Ehrhart for polygons}
L_\P(t) 
&=
{\rm area }(\P) t^2 
-
\frac{1}{2} 
\sum_{v \in V}  
\left(
\frac{\langle v, z \rangle }{\langle w_1(v), z \rangle}
+\frac{\langle v, z \rangle }{\langle w_2(v), z \rangle}
\right) t
+1,
\end{align}
valid for almost all $z \in \C^2$ and 
$t \in \Z_{>0}$.
To show how this follows from Corollary 
\ref{thm:Ehrhart.from.Bernoulli.Barnes},
we first recall the definition 
$\va_v:=  \icol{   \langle  w_1(v), z \rangle \\   \langle w_2(v), z \rangle }$ at each vertex $v$ of $\P$.  By Corollary 
\ref{thm:Ehrhart.from.Bernoulli.Barnes}, we have:
\begin{align*}
L_\P(t)  
&=     
\frac{1}{2} 
\sum_{k=0}^2  
\binom{2}{k}
\sum_{v \in V}                             
B_k(t  \langle v, z \rangle,   \va_{v})
   \sum_{p \in \Pi_v \cap \Z^2}   
   {\langle  p, z \rangle}^{2-k}\\
&=
\frac{1}{2} 
\sum_{v \in V}                             
B_0(t  \langle v, z \rangle,   \va_{v})
   \sum_{p \in \Pi_v \cap \Z^2}   
   {\langle  p, z \rangle}^{2} +
\sum_{v \in V}                             
B_1(t  \langle v, z \rangle,   \va_{v})
   \sum_{p \in \Pi_v \cap \Z^2}   
   {\langle  p, z \rangle}\\
&+\frac{1}{2}\sum_{v \in V}                      B_2(t  \langle v, z \rangle,   \va_{v})
   \vol \Pi_v\\
&=
\left(
\frac{1}{2}\sum_{v \in V}
\frac{\sum_{p \in \Pi_v \cap \Z^2}{\langle p,z\rangle}^{2}}
{\langle w_1(v),z\rangle \langle w_2(v),z\rangle}
-\frac{1}{2}\sum_{v \in V}
\frac{\langle w_1(v),z\rangle+\langle w_2(v),z\rangle}
{\langle w_1(v),z\rangle \langle w_2(v),z\rangle}
\sum_{p \in \Pi_v \cap \Z^2}\langle p,z\rangle
\right.\\
&\left.\qquad
+\frac{1}{12}\sum_{v \in V}\vol \Pi_v
\left(
\frac{\langle w_1(v),z\rangle}{\langle w_2(v),z\rangle}
+
\frac{\langle w_2(v),z\rangle}{\langle w_1(v),z\rangle}
\right)
+\frac{1}{4}\sum_{v \in V}\vol \Pi_v
\right)\\
&\quad
+t\left(
\sum_{v \in V}
\frac{\langle v,z\rangle}
{\langle w_1(v),z\rangle \langle w_2(v),z\rangle}
\sum_{p \in \Pi_v \cap \Z^2}\langle p,z\rangle
-\frac{1}{2}\sum_{v \in V}
\langle v,z\rangle \vol \Pi_v
\left(
\frac{1}{\langle w_1(v),z\rangle}
+
\frac{1}{\langle w_2(v),z\rangle}
\right)
\right)\\
&\quad
+t^2\left(
\frac{1}{2}\sum_{v \in V}
\langle v,z\rangle^2
\frac{\vol \Pi_v}
{\langle w_1(v),z\rangle \langle w_2(v),z\rangle}
\right),
\end{align*}
where we've used \eqref{ex: dim 2, B_0, B_1} 
and \eqref{ex: dim 2, B_2} for the three Barnes polynomials above. 
To simplify notation, we define 
 \[
 c_v:= \frac{1}{2}(w_1(v) + w_2(v)),
 \]
the center of mass of the continuous vertex parallelepiped $\Pi_v$. 
To summarize:
\begin{align}
\label{pre-summary of L_P for polygon}
L_\P(t) &= 
{\rm area }(\P) t^2 + 
t
\sum_{v \in V}  
\frac{\langle v, z \rangle}  
{\langle w_1(v), z \rangle \langle w_2(v), z \rangle}
\left(
{\displaystyle\sum_{p \in \Pi_v \cap \Z^2} {\langle  p, z \rangle}}
-\vol \Pi_v \langle c_v, z \rangle
\right)+1,
\end{align}
valid for almost all $z\in \C^2$ and 
$t \in \Z_{>0}$.
We used Corollary \ref{cor:ConstantCoeff} 
for the constant term in \eqref{pre-summary of L_P for polygon},
and equation  \eqref{Brion volume} for the area of the polygon.   
It's an easy exercise that summing all the points in a $2$-dimensional vertex parallelepiped gives us:
\begin{equation}
\sum_{p \in \Pi_v \cap \Z^2}  p
= \big(\vol \Pi_v -1  \big)\ c_v.
\end{equation}
[hint: translate $\Pi_v \cap \Z^2$ to the origin by (its center of mass) $c_v$, and notice that the new non-zero points cancel in pairs].
So we have:
\[
\sum_{p \in \Pi_v \cap \Z^2} 
\langle  p, z \rangle - \vol \Pi_v\langle c_v, z\rangle
= -\langle c_v, z\rangle.
\]
Substituting the latter expression into 
\eqref{pre-summary of L_P for polygon}, we obtain
\begin{align}
L_\P(t) 
&= 
{\rm area }(\P) t^2 
- 
t \sum_{v \in V}  
\frac{\langle v, z \rangle \langle c_v, z\rangle}  
{\langle w_1(v), z \rangle \langle w_2(v), z \rangle}
+1 \\
&= \label{Coordinate free polygon Ehrhart}
{\rm area }(\P) t^2 
-
\frac{1}{2} 
\sum_{v \in V}  
\left(
\frac{\langle v, z \rangle }{\langle w_1(v), z \rangle}
+\frac{\langle v, z \rangle }{\langle w_2(v), z \rangle}
\right) t
+1,
\end{align}
using the definition $\langle c_v, z \rangle 
= \frac{1}{2}\left( \langle w_1(v),z \rangle +\langle w_2(v),z \rangle\right)$.

Finally, let's see how to translate the formulation 
\eqref{Coordinate free polygon Ehrhart} into the standard known result on Ehrhart polynomials.  
For any $x\in \Z^d$, we define $\gcd(x):= |\gcd(x_1, \dots, x_d)|$.   
We consider any four consecutive vertices of $\P$, say $v_1, v_2, v_3, v_4$.  
To give the primitive edge vectors for $v_2$, we define 
$g_{21}:= \gcd(v_1 - v_2)$, which is the number of integer points on the half-open edge joining $v_1$ and $v_2$.  We note the important fact that  $g_{23} = g_{32}$. 
Then $w_1(v_2)= \frac{1}{g_{21}} (v_1 - v_2)$, and similarly
$w_2(v_2)= \frac{1}{g_{23}} (v_3 - v_2)$.  In the summation of \eqref{Coordinate free polygon Ehrhart}, we focus on just the two vertices $v_2$ and $v_3$, and we consider their contribution to the sum:
\begin{align}
&\left(
\frac{\langle v_2, z \rangle }{\langle w_1(v_2), z \rangle}
+\frac{\langle v_2, z \rangle }{\langle w_2(v_2), z \rangle}
\right)
+
\left(
\frac{\langle v_3, z \rangle }{\langle w_1(v_3), z \rangle}
+\frac{\langle v_3, z \rangle }{\langle w_2(v_3), z \rangle}
\right)\\ 
\label{sum two polygon vertices}
&=\frac{g_{21}\langle v_2, z \rangle }{\langle v_1 - v_2, z \rangle}
+
\left(
\frac{g_{23}\langle v_2, z \rangle }{\langle v_3 - v_2, z \rangle}
+\frac{g_{32}\langle v_3, z \rangle }{\langle v_2 - v_3, z \rangle}
\right)
+\frac{g_{34}\langle v_3, z \rangle }{\langle v_4 - v_3, z \rangle}\\ 
\label{simplify the sum of two polygon vertices} 
&= 
\frac{g_{21}\langle v_2, z \rangle }{\langle v_1 - v_2, z \rangle}
-g_{23}
+\frac{g_{34}\langle v_3, z \rangle }{\langle v_4 - v_3, z \rangle}.  
\end{align}
Once we consider the full summation in \eqref{Coordinate free polygon Ehrhart}, we therefore see that this simplification in pairs, going from \eqref{sum two polygon vertices} 
to \eqref{simplify the sum of two polygon vertices},  gives us the classical result:
\begin{align}
L_\P(t) 
&= 
{\rm area }(\P) t^2 
+\frac{1}{2}
t \sum_{j=1}^N \gcd \left(v_j - v_{j+1}\right)
+1,
\end{align}
where $N$ is the number of vertices of $\P$, and by definition $v_{N+1} := v_1$. 


\section{Example: a rational triangle}
\label{Example:rational triangle}
\begin{figure}[hbt!]
    \centering
\includegraphics[scale=0.5]{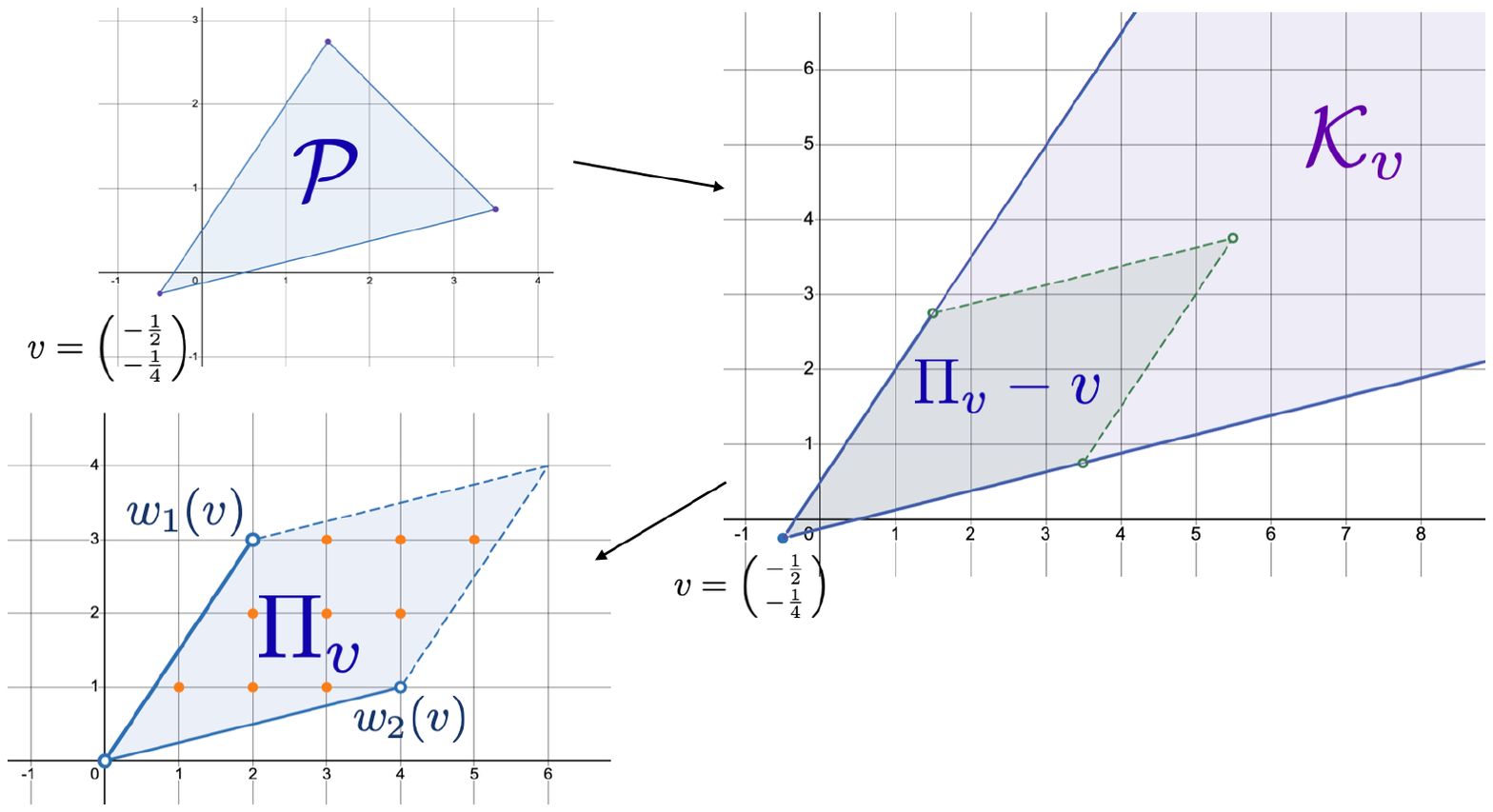}
\caption{
{\bf Top left}: a rational triangle $\P$ with a rational vertex $v$. 
{\bf right}: the vertex tangent cone $\K_v$ at the vertex $v$. 
{\bf Bottom left}: the vertex parallelepiped $\Pi_v$, with its integer edge vectors $w_1(v), w_2(v)$.  Note that $\Pi_v$ always has a vertex at the origin. This picture is  standard in the geometry of numbers (see for example \cite{Sinai2})
}
\label{Fig:triangle wit Pi_v}
\end{figure}

\begin{figure}[hbt!]
    \centering
\includegraphics[scale=0.55]{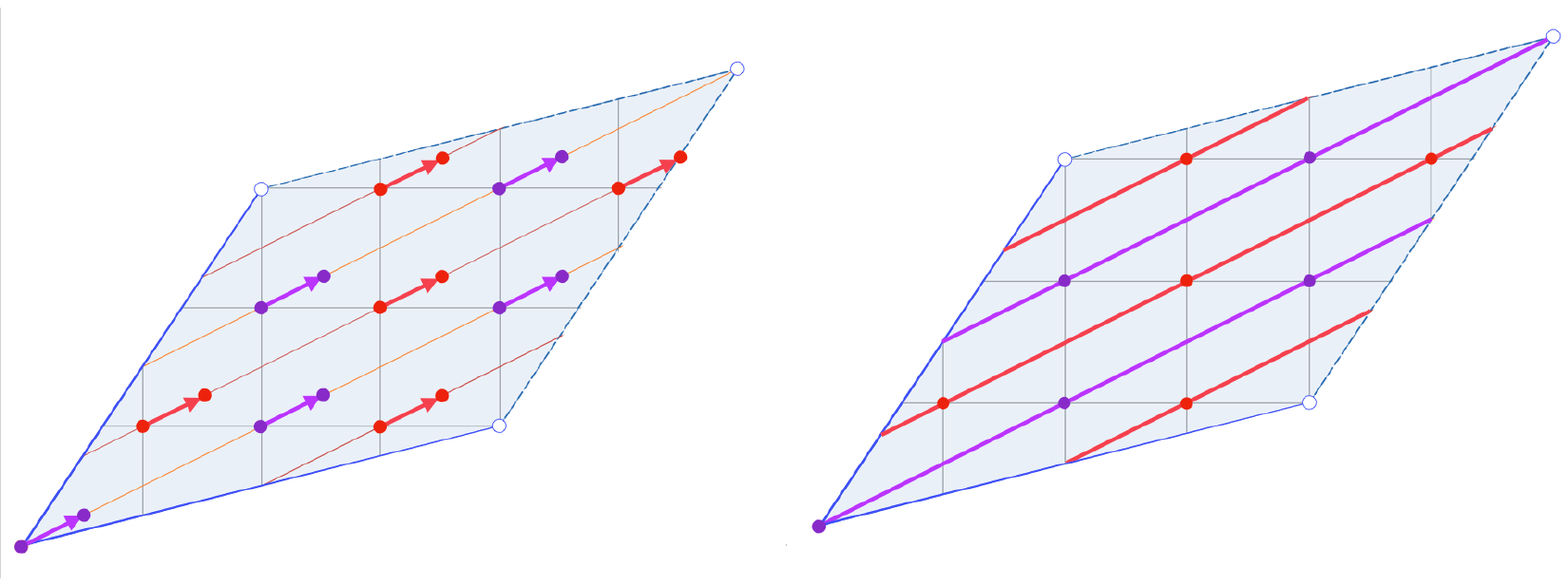}
\caption{{\bf Left}:  The vertex parallelepiped $\Pi_v$ for the vertex $v_0:= \icol{-\frac{1}{2} \\ -\frac{1}{4} }$.  The lattice flow 
$\left\{ \left(\Z^d - tv_0 \right) \mod \Pi_{v_0} \ \mid \ 0 \leq t \leq \frac{1}{4} 
\right\}$ is drawn on the left with red and purple line segments. 
Each integer point is flowing along its geodesic on this compact $2$-torus, from $t=0$ until $t = \frac{1}{4}$.
 {\bf Right}: the same lattice flow has reached its full closed geodesic flow, at $t=4$. 
    }
    \label{Fig:lattice flow 1}
\end{figure}

We consider the rational triangle $\P$ drawn in Figure 
\ref{Fig:triangle wit Pi_v}, whose three vertices are given by 
$v_0:= 
\icol{-\frac{1}{2} \\ -\frac{1}{4} }, 
v_1:= 
\icol{ \frac{7}{2}\\  \frac{3}{4}},   
v_2:= \icol{ \frac{3}{2} \\ \frac{11}{4}}$. 
The partial lattice flow shown on the left of Figure \ref{Fig:lattice flow 1}
is the union of the purple and red line segments,  and is  defined by 
\[
\left\{ 
LatticeFlow(t) \ \mid \ 0 \leq t \leq \frac{1}{4} 
\right\}:=
\left\{ 
\Z^2 - tv \pmod {\Pi_v} \ \big| \ 
0\leq t \leq \frac{1}{4}.
\right\}
\]
Each integer point is flowing along its geodesic defined on the flat (and compact) $2$-torus given by $\Pi_v$.
The lattice flow reaches its full closed geodesic flow, at $t=4$, as shown on the right-hand side of Figure
\ref{Fig:lattice flow 1}.
We note that there are exactly two closed geodesics here (one in red and the other in purple). Considering the integer points in $\Pi_v$ as a  finite abelian group $G$ under addition mod $\Pi$, we see that the origin is contained in a subgroup of $G$ that has order $5$, and all closed geodesics arise from the cosets of this subgroup.

\noindent
We also observe that for any lattice flow, the orbit of any integer point in $\Pi_v$ produces a closed geodesic that has the same length as any other integer point.  

In Figure \ref{Fig:lattice flow 1}, there are exactly two closed geodesics, and the total length of the red geodesic equals the total length of the blue geodesic, because here the blue geodesic arises from an index $2$ subgroup of $\Z^2/\Lambda$.


\noindent
{\bf The vertex $v_0$}. Its two primitive edge vectors are
$w_1(v_0)=\icol{2\\3}$, and 
$w_2(v_0)=\icol{4\\1}$. 
We have:
\begin{equation}
W_{v_0}
:=
\begin{pmatrix}
2&4\\
3&1
\end{pmatrix},
\qquad
\Lambda_{v_0}
:=
W_{v_0}\Z^2,
\qquad
\Pi_{v_0}
:=
W_{v_0}[0,1)^2,
\label{eq:rational-triangle-v0-matrix-lattice}
\end{equation}
with $|\det W_0|= \vol(\Pi_{v_0})=10$.
The ten integer points in $\Pi_{v_0}$ are:
\begin{equation}
\mathcal R_{v_0}
=
\bigl\{
(0,0),(1,1),(2,1),(2,2),(3,1),
(3,2),(3,3),(4,2),(4,3),(5,3)
\bigr\}.
\end{equation}
These points represent the finite abelian group
$G_{v_0}
:=
\Z^2/\Lambda_{v_0}$, which in this case is a cyclic group. 
The direction of the lattice flow is
$-v_0
=
\icol{\frac{1}{2}\\\frac{1}{4}}$. 
After time $t=4$, the displacement is the integer vector
$
-4v_0
=
\icol{2\\1}$.
Consequently, the equivalence class
of $h:=
\left[
\icol{2\\1}
\right]$
has order $5$.  It generates the index $2$ subgroup
\begin{equation}
H
:=
\langle h\rangle
=
\bigl\{
[(0,0)],[(2,1)],[(4,2)],
[(2,2)],[(4,3)]
\bigr\}.
\label{eq:rational-triangle-order-five-subgroup}
\end{equation}
The other coset is
$
g+H
=
\bigl\{
[(1,1)],[(3,2)],[(5,3)],
[(3,3)],[(3,1)]
\bigr\}$. 
So  $G_{v_0}=
H\sqcup(g+H)$.
It's interesting to describe the flow explicitly.  
Since
$W_{v_0}^{-1}(-v_0)
=
\icol{\frac{1}{20}\\\frac{1}{10}}$,
a particular point $q \in \Pi_{v_0} \cap \Z^2$ has the flow trajectory 
\begin{equation}
q-tv_0  \  (\bmod\Pi_v)
=
W_{v_0}
\left\{
W_{v_0}^{-1}q
+
t\icol{\frac{1}{20}\\\frac{1}{10}}
\right\}.
\label{eq:rational-triangle-explicit-trajectory}
\end{equation}

We next evaluate Theorem~1 at all three vertices for positive integer dilations. We use
$n\in\mathbb Z_{>0}$ for the integer dilation parameter, and make the convenient
global choice

\begin{equation}
z:=e_1=(1,0).
\label{eq:triangle-choice-z}
\end{equation}

For this choice of $z$, the primitive edge data at the three vertices are

\begin{equation}
\begin{array}{c|c|c|c|c}
j
&
w_1(v_j)
&
w_2(v_j)
&
a_j
&
\langle v_j,e_1\rangle
\\ \hline
0
&
(2,3)
&
(4,1)
&
(2,4)
&
-\dfrac12
\\[1mm]
1
&
(-4,-1)
&
(-1,1)
&
(-4,-1)
&
\dfrac72
\\[1mm]
2
&
(-2,-3)
&
(1,-1)
&
(-2,1)
&
\dfrac32
\end{array}
\label{eq:triangle-three-local-data}
\end{equation}

The corresponding vertex parallelepiped volumes are
$\vol(\Pi_0)=10,
\vol(\Pi_1)=5$, and
$\vol(\Pi_2)=5$.
We now compute the shifted polytope Dedekind sums that occur in Theorem~1.
Since
$4v_0, 4v_1, 4v_2\in\mathbb Z^2$, 
each of the three lattice flows has period $4$ as an unlabelled set.
Consequently,
\begin{equation}
\mu_k(\Pi_j,e_1,(n+4)v_j)
=
\mu_k(\Pi_j,e_1,nv_j),
\qquad
j\in\{0,1,2\}.
\label{eq:triangle-moments-period-four}
\end{equation}

Thus it is enough to compute the discrete moments for the four residue classes
of $n$ modulo $4$. Writing

\begin{equation}
r\equiv n\pmod 4,
\qquad
0\leq r\leq3,
\label{eq:triangle-residue-r}
\end{equation}
direct reduction of the shifted lattice points into their respective half-open
vertex parallelepipeds gives
\begin{equation}
\resizebox{\textwidth}{!}{$
\begin{array}{c|cc|cc|cc}
r
&
\mu_1(\Pi_0,e_1,nv_0)
&
\mu_2(\Pi_0,e_1,nv_0)
&
\mu_1(\Pi_1,e_1,nv_1)
&
\mu_2(\Pi_1,e_1,nv_1)
&
\mu_1(\Pi_2,e_1,nv_2)
&
\mu_2(\Pi_2,e_1,nv_2)
\\ \hline
0
&27&93
&-10&30
&-2&2
\\[1mm]
1
&28&\dfrac{189}{2}
&-\dfrac{23}{2}&\dfrac{133}{4}
&-\dfrac52&\dfrac{13}{4}
\\[2mm]
2
&27&89
&-12&34
&-3&3
\\[1mm]
3
&28&\dfrac{197}{2}
&-\dfrac{27}{2}&\dfrac{173}{4}
&-\dfrac72&\dfrac{21}{4}
\end{array}
$}
\label{eq:triangle-moment-table}
\end{equation}

The zeroth polytope Dedekind sums are independent of $n$, and are simply the
numbers of lattice points in the corresponding vertex parallelepipeds:

\begin{equation}
\mu_0(\Pi_0,e_1,nv_0)=10,
\qquad
\mu_0(\Pi_1,e_1,nv_1)=5,
\qquad
\mu_0(\Pi_2,e_1,nv_2)=5.
\label{eq:triangle-zero-moments}
\end{equation}

For each $j\in\{0,1,2\}$, Theorem~1 gives the vertex contribution
\begin{align}
C_j(n)
&:=
\frac12\Big(
B_0\!\left(
n\langle v_j,e_1\rangle,
a_j
\right)
\mu_2(\Pi_j,e_1,nv_j)
+
2B_1\!\left(
n\langle v_j,e_1\rangle,
a_j
\right)
\mu_1(\Pi_j,e_1,nv_j)
\\
&+ 
B_2\!\left(
n\langle v_j,e_1\rangle,
a_j
\right)
\mu_0(\Pi_j,e_1,nv_j)
\Big).
\label{eq:triangle-vertex-contribution}
\end{align}

For example, suppose that $n\equiv0\pmod4$. At the vertex $v_0$, we have
$\mu_1(\Pi_0,e_1,nv_0)=27, \ 
\mu_2(\Pi_0,e_1,nv_0)=93$, and
$\mu_0(\Pi_0,e_1,nv_0)=10$.
Moreover, since $a_0=(2,4)$ and
$\langle v_0,e_1\rangle=-1/2$, the explicit Barnes polynomials in
(68) and (69) give

\begin{equation}
B_0\!\left(-\frac n2,(2,4)\right)
=
\frac18,
\label{eq:v0-B0}
\end{equation}

\begin{equation}
B_1\!\left(-\frac n2,(2,4)\right)
=
-\frac{n}{16}-\frac38,
\label{eq:v0-B1}
\end{equation}

and

\begin{equation}
B_2\!\left(-\frac n2,(2,4)\right)
=
\frac{n^2}{32}
+
\frac{3n}{8}
+
\frac{11}{12}.
\label{eq:v0-B2}
\end{equation}

Therefore

\begin{equation}
\begin{aligned}
C_0(n)
&=
\frac12
\left[
\frac{93}{8}
+
54
\left(
-\frac{n}{16}-\frac38
\right)
+
10
\left(
\frac{n^2}{32}
+\frac{3n}{8}
+\frac{11}{12}
\right)
\right]
\\
&=
\frac{15n^2+18n+26}{96}.
\end{aligned}
\label{eq:triangle-C0-residue-zero}
\end{equation}

Thus, at the vertex $v_0$, we have seen explicitly how the lattice flow in
Figure~2 produces the first and second polytope Dedekind sums, and how these
in turn enter the Barnes-polynomial formula of Theorem~1.

Performing the corresponding substitutions in
\eqref{eq:triangle-vertex-contribution} at all three vertices gives the
following vertex contributions:

\begin{equation}
\begin{array}{c|ccc}
r
&
96C_0(n)
&
96C_1(n)
&
96C_2(n)
\\ \hline
0
&
15n^2+18n+26
&
735n^2+210n+50
&
-270n^2-36n+20
\\
1
&
15n^2+12n-1
&
735n^2+84n-1
&
-270n^2+2
\\
2
&
15n^2+18n+2
&
735n^2+42n-22
&
-270n^2+36n+20
\\
3
&
15n^2+12n+23
&
735n^2-84n-1
&
-270n^2+72n-22
\end{array}
\label{eq:triangle-contribution-table}
\end{equation}

Although the individual vertex contributions arise from rational expressions
in the auxiliary $z$-parameter, their sum is independent of $z$, as
Theorem~1 guarantees. For our choice $z=e_1$, we therefore have

\begin{equation}
L_P(n)
=
C_0(n)+C_1(n)+C_2(n).
\label{eq:triangle-sum-contributions}
\end{equation}
The area of the triangle is
$5$, which is the coefficient of $n^2$. Summing the three
vertex contributions in \eqref{eq:triangle-contribution-table} gives us:
\begin{equation}
L_P(n)
=
\begin{cases}
5n^2+2n+1,
&
n\equiv0\pmod4,
\\[1mm]
5n^2+n,
&
n\equiv1\pmod4,
\\[1mm]
5n^2+n,
&
n\equiv2\pmod4,
\\[1mm]
5n^2,
&
n\equiv3\pmod4.
\end{cases}
\label{eq:triangle-final-ehrhart-quasipolynomial}
\end{equation}

We have therefore recovered the Ehrhart quasi-polynomial of $P$ directly
from the lattice flows in its three vertex parallelepipeds. Its period divides
$4$, since $4v_j\in\mathbb Z^2$ for every vertex. On the other hand, the
constituent polynomials in 
\eqref{eq:triangle-final-ehrhart-quasipolynomial} show that the Ehrhart
quasi-polynomial of $P$ has exact period $4$.


\section{Example: a d-dimensional cube}
\label{The cube}
We consider the example of a 
$d$-dimensional cube $\square:= [-1, 1]^d$, another smooth polytope.   Although it is
 trivial to compute the  Ehrhart polynomial in this case, namely $L_\square(t) = (2t+1)^d$, the verification below will shed some light on our use of the Barnes polynomials, and even gives new identities
 (equation \eqref{identity for the d-cube} below)  for the discrete hypercube $\{ -1, 1\}^d$.

To begin, we label each of the $2^d$ vertices  
$v \in \square$ by associating to each vertex a choice of signs as follows:  
\begin{equation}
    v :=  \icol{ \alpha_1 \\  . \\ . \\ . \\ \alpha_d},
\end{equation}
where we fix each $\alpha_k  \in \{-1, 1\}$.  It's clear that for any such fixed vertex $v$, the  primitive integer edge vectors that are incident with $v$ are:
 $w_1(v) = - \alpha_1 \e_1,  \dots,  w_d(v) =   - \alpha_d \e_d$, where we are using the standard basis vectors $\e_1, \dots, \e_d$.  It follows that the vector of linear forms associated with the vertex $v$ is:
 \begin{align*}
\va_v &:= \Big( \langle  w_1(v), z \rangle, \cdots, \langle w_d(v), z \rangle \Big)^T\\
&= \left(   \langle  - \alpha_1 \e_1, z \rangle, \cdots,   \langle  - \alpha_d \e_d, z \rangle \right)^T \\
&= - \left(   \alpha_1 z_1, \cdots,      \alpha_d z_d     \right)^T. 
 \end{align*}
By Corollary \ref{Ehrhart for smooth polytopes}, we have
\begin{align*}
L_{\square}(t)   &=  \frac{  (-1)^d }{d!}   \sum_{j=0}^d    \sum_{v \in \rm{V}}
   \left(\vol \Pi_{v}\right)     
\binom{d}{j}
 B_{d-j}(\va_{v}) 
       \langle v, z \rangle^{j}  \  t^{j} \\
&=  \frac{  (-1)^d }{d!}    
\sum_{j=0}^d    
\sum_{(\alpha_1, \cdots, \alpha_d) \in \{ -1, 1 \}^d}
B_{d-j}(\va_{v})   
\Big(
\alpha_1 z_1 + \cdots + \alpha_d z_d
\Big)^{j}   
\binom{d}{j}  t^{j},
\end{align*}
where we are summing over the finite hypercube $\{ -1, 1\}^d$, consisting of $2^d$ elements.
Comparing this result with the known Ehrhart polynomial
$L_{\square}(t) = (2t+1)^d = 
\sum_{j=0}^d 2^j t^j \binom{d}{j}$, we observe the following
identities for the Barnes numbers, and for the discrete hypercube:
\begin{equation} \label{identity for the d-cube}
  \frac{  (-1)^d }{d!}    \sum_{(\alpha_1, \cdots, \alpha_d) \in \{ -1, 1 \}^d}
\Big(
\alpha_1 z_1 + \cdots + \alpha_d z_d
\Big)^{j}   B_{d-j} \icol{ -\alpha_1 z_1 \\ .\\ .\\ .\\  -\alpha_d z_d  }   = 2^j,
 \end{equation}
for each $0\leq j\leq d$, and for almost all $z \in \C^d$.  

We may notice that the case $j=0$ of
identity \eqref{identity for the d-cube}
was already covered by  Corollary \ref{cor:ConstantCoeff}, in the special case of smooth polytopes.  The case $j=d$ of identity \eqref{identity for the d-cube}
seems interesting, and using 
$B_{0}(\va_{v}) =    \frac{1}{     \prod_{k=1}^d  \langle w_k(v), z \rangle }$, tells us that 
\begin{align}
2^d &=  \frac{  (-1)^d }{d!}    \sum_{(\alpha_1, \cdots, \alpha_d) \in \{ -1, 1 \}^d}
      (\alpha_1 z_1 + \cdots + \alpha_d z_d)^{d}   B_{0}(\va_{v})  \\
      &=   \frac{  (-1)^d }{d!}    \sum_{(\alpha_1, \cdots, \alpha_d) \in \{ -1, 1 \}^d}
      \frac{    (\alpha_1 z_1 + \cdots + \alpha_d z_d)^{d}      }{\prod_{k=1}^d \langle w_k(v), z \rangle} \\
      &= \frac{  (-1)^d }{d!}    \sum_{(\alpha_1, \cdots, \alpha_d) \in \{ -1, 1 \}^d}
      \frac{    (\alpha_1 z_1 + \cdots + \alpha_d z_d)^{d}      }{\prod_{k=1}^d (-\alpha_k z_k) }  \\
      &= \frac{  1}{d!}   \sum_{(\alpha_1, \cdots, \alpha_d) \in \{ -1, 1 \}^d}
      \frac{    (\alpha_1 z_1 + \cdots + \alpha_d z_d)^{d}      }{\prod_{k=1}^d \alpha_k z_k },
\end{align}
for almost all $z\in \C^d$.
As we've seen above,  we can also reverse-engineer our approach.  Namely,  
 if we have an integer polytope $\P$ for which we already know
$L_\P(t)$, then the polytope $\P$ determines a set of identities for the Barnes numbers, as we saw in \eqref{identity for the d-cube}.

\begin{example}
In $\R^2$, just for fun let's retrieve the Ehrhart polynomial for the standard triangle $\P:= \conv \{ \icol{ 0 \\ 0}, \icol{ 1 \\ 0},   \icol{ 0 \\ 1}   \}$, using Barnes polynomials. We'll use  Corollary \ref{Ehrhart for smooth polytopes} of Theorem
\ref{thm:Ehrhart.from.Bernoulli.Barnes}, 
since the standard triangle is a smooth polygon. 
Here we only need to use $B_2(x, \va)$, the second Barnes polynomial, given by \eqref{ex: dim 2, B_2}.  By definition we have 
\[
\va_{v_1}:=  
\icol{  \langle z,  \   \omega_1(v_1) \rangle 
 \\ \langle z, \    \omega_2(v_1) \rangle}
=  \icol{ z_1 \\ z_2}, 
\quad
\va_{v_2}:= 
\icol{   \langle z, \  \omega_1(v_2) \rangle
 \\ \langle z, \    \omega_2(v_2) \rangle}
= \icol{ -z_1 + z_2\\ -z_1},
\] 
and 
$\va_{v_3}:= 
 \icol{ z_1 - z_2\\ -z_2}$.
By Corollary \ref{Ehrhart for smooth polytopes}
we have
\begin{align*}
&  2! L_P(t) =(-1)^2 \Big(   B_2(t\cdot 0, \va_{v_1}) + B_2(   tz_1,  \va_{v_2} )
+  B_2(   tz_2,  \va_{v_3} )   \Big) \\
&=  \frac{1}{6} \left( \tfrac{z_1}{z_2} + \tfrac{z_2}{z_1} \right) + \tfrac{1}{2}\\
&+ t^2 z_1^2 \left( \tfrac{1}{(-z_1+z_2)(-z_1)} \right) 
-tz_1 \left( \tfrac{1}{-z_1 + z_2} + \tfrac{1}{-z_1} \right)
+ \tfrac{1}{6} \left(  \tfrac{-z_1 + z_2}{-z_1} +  \tfrac{-z_1}{-z_1 + z_2} \right) + \tfrac{1}{2} \\
&+ t^2 z_2^2 \left( \tfrac{1}{(-z_2+z_1)(-z_2)} \right) 
-tz_2 \left( \tfrac{1}{-z_2 + z_1} + \tfrac{1}{-z_2} \right)
+ \tfrac{1}{6} \left(  \tfrac{-z_2 + z_1}{-z_2} +  \tfrac{-z_2}{-z_2 + z_1} \right) + \tfrac{1}{2} \\
&= t^2 + 3t + 2,
\end{align*}
and we see that all of the $z$'s miraculously cancel out, as predicted by Theorem  \ref{thm:Ehrhart.from.Bernoulli.Barnes}.
\hfill $\hexago$
\end{example}

\section{Proof of Corollary
\ref{Cor:novel formula for codim 1 coeff}}

\begin{proof}
\hypertarget{Codim 1 formula}
{(of Corollary \ref{Cor:novel formula for codim 1 coeff})}
For each vertex $v$, we recall the edge matrix $W_v$ defined in \eqref{def:edge matrix}.
With  $q_{v,k}$ and $\beta_{v,k}$  as in the statement of the
corollary, we'll also use the first two Barnes numbers,
\[
B_0(\va)=\frac{1}{a_1\cdots a_d},
\qquad
B_1(\va)
=-\frac12\frac{a_1+\cdots+a_d}{a_1\cdots a_d}.
\]
By Theorem~\ref{thm:Ehrhart quasi-polynomial from Barnes}, part
\eqref{part b of Theorem 2}, with $r=d-1$, we have
\begin{align}
c_{d-1}(t)
&=
\frac{(-1)^d}{(d-1)!}
\sum_{v\in V}
\langle v,z\rangle^{d-1}
\left(
B_0(\va_v)
\sum_{p\in\Z^d-tv\;(\bmod\Pi_v)}
\langle p,z\rangle
+
B_1(\va_v)\vol(\Pi_v)
\right)\notag\\
&=
\frac{(-1)^d}{(d-1)!}
\sum_{v\in V}
\frac{\langle v,z\rangle^{d-1}}
{\displaystyle\prod_{j=1}^d\langle w_j(v),z\rangle}
\left(
\sum_{p\in\Z^d-tv\;(\bmod\Pi_v)}
\langle p,z\rangle
-
\frac12\vol(\Pi_v)
\sum_{k=1}^d\langle w_k(v),z\rangle
\right).
\label{eq:codim-one-before-first-moment}
\end{align}
Lemma~\ref{lem:closed form, first discrete moment} gives
\begin{align*}
&\sum_{p\in\Z^d-tv\;(\bmod\Pi_v)}
\langle p,z\rangle
-
\frac12\vol(\Pi_v)
\sum_{k=1}^d\langle w_k(v),z\rangle\\
&\quad=
\vol(\Pi_v)
\sum_{k=1}^d
\left[
\frac{1}{q_{v,k}}
\left(
\frac{q_{v,k}-1}{2}
+
\left\{-tq_{v,k}\beta_{v,k}\right\}
\right)
-
\frac12
\right]
\langle w_k(v),z\rangle\\
&\quad=
-\vol(\Pi_v)
\sum_{k=1}^d
\frac{1}{q_{v,k}}
\left(
\frac12-
\left\{-tq_{v,k}\beta_{v,k}\right\}
\right)
\langle w_k(v),z\rangle.
\end{align*}
Substituting this identity into
\eqref{eq:codim-one-before-first-moment} and cancelling the factor
$\langle w_k(v),z\rangle$ gives
\begin{align*}
c_{d-1}(t)
&=
\frac{(-1)^{d-1}}{(d-1)!}
\sum_{v\in V}
\sum_{k=1}^d
\frac{
\langle v,z\rangle^{d-1}\vol(\Pi_v)
}{
\displaystyle\prod_{\substack{1\leq j\leq d\\j\neq k}}
\langle w_j(v),z\rangle
}
\frac{1}{q_{v,k}}
\left(
\frac12-
\left\{-tq_{v,k}\beta_{v,k}\right\}
\right).
\end{align*}
\end{proof}


\section{Proof of Theorem  
\ref{thm:Ehrhart quasi-polynomial from Barnes}}
When a vertex $v\in \P$ is an integer point, the identity
$\sigma_{ \Pi_v+ v}(z)  = e^{\langle v, z \rangle} \sigma_{\Pi_v}(z)$ holds.
When $v \notin \Z^d$, however,  the latter identity is no longer necessarily true.
Therefore more care has  to be taken in case we consider the vertex parallelepiped of a rational vertex, or in case we replace $v$ with $tv$, for $t \in \Q$.

We also recall the standard fact that the integer point enumerator for the parallelepiped 
$\Pi_v + v$ is defined by
\begin{equation}
\sigma_{\Pi_v+ v}(z) := \sum_{p \in  ( \Pi_v + v )   \cap \Z^d}  e^{\langle p, z \rangle},
\end{equation}
for all $z \in \C^d$.
Given a  $d$-dimensional simplicial cone $\K_v \subset \R^d$,  with any real apex $v \in \R^d$, and with $d$
 linearly independent integer edge vectors $\omega_1, \dots, \omega_d \in \Z^d$, we have:
\begin{equation} \label{integer point transform of a cone}
\sum_{n \in \K_v \cap \Z^d}    e^{ \langle n,  z \rangle}
= 
\frac{   \sigma_{ \Pi_v + v}(z)     }
{       \prod_{k=1}^d  \left(      1 - e^{ \langle  \omega_k   , z \rangle}   \right)          }.
\end{equation}
\cite[Theorem 3.4]{BeckRobins}.
We use the notation $\K_v - v := \K_0(v)$, the translation of the cone $\K_v$ to the origin.  We observe that $t\K_0(v) = \K_0(v)$ for all $t>0$, a trivial but useful fact.  Consequently, when we dilate
$\K_v$ by any $t>0$, \emph{all of its edge vectors remain the same}, and only its apex $v$ gets dilated:  
\begin{equation}
\label{trivial dilation of a cone}
t\K_v = t(\K_v - v) + tv = 
t\K_0(v) + tv = \K_0(v) + tv.  
\end{equation}
We  denote the latter cone by $\K_{tv}:= \K_0(v) +tv$.

\begin{proof}
\hypertarget{Proof of Theorem 2}
{(of Theorem \ref{thm:Ehrhart quasi-polynomial from Barnes})}
We begin with the discrete Brion Theorem:
\begin{equation}\label{first general step}
\sigma_{t\P}(z) = \sigma_{t\K_{v_1}}(z)  + \cdots +  \sigma_{t\K_{v_N}}(z),
\end{equation}
an identity that is valid for almost all $z\in \C^d$, and valid for any rational $t>0$. 
For Steps $1$--$3$ below, we first fix $t\in\Q_{>0}$.  After deriving the desired formula for all positive rational values of $t$, we extend it to every real $t>0$ by a piecewise-polynomial argument.  We will also use the observation \eqref{trivial dilation of a cone} regarding the meaning of $\K_{tv}$ for all positive $t>0$.

\noindent
{\bf Step $1$}.  
From  \eqref{trivial dilation of a cone}, it follows that the elementary identity \eqref{integer point transform of a cone} becomes:
\begin{equation} \label{integer point transform of a real dilated cone}
\sigma_{t\K_{v}}(z):=
\sum_{n \in t\K_v \cap \Z^d}    
e^{ \langle n,  z \rangle}
= 
\frac{   \sigma_{ \Pi_v + tv}(z) }
{       \prod_{k=1}^d  \left(      1 - e^{ \langle  \omega_k   , z \rangle}   \right)          },
\end{equation} 
for all $t>0$.
We now substitute $z:= x z_0$, with $x>0$ and $z_0$ a generic direction (thus avoiding any singularities in the denominators), and we begin to expand $\sigma_{t\K_{v}}(xz_0)$ as a Taylor series about $x=0$:
\begin{align}
\sigma_{t\K_{v}}(xz_0)  
&= 
\frac{ \sigma_{ \Pi_v + tv}(xz_0) }
{  \prod_{k=1}^d \Big( 1 - e^{x\langle z_0, w_k(v) \rangle}  \Big) } \\
&= 
\frac{ e^{-t x\langle v, z_0 \rangle} \sigma_{ \Pi_v + tv}(xz_0) }{x^d}
\frac{ x^d e^{t x\langle v, z_0 \rangle} }
{  \prod_{k=1}^d \Big( 1 - e^{x\langle z_0, w_k(v) \rangle}  \Big) } \\
&=
\frac{  e^{-t x\langle v, z_0 \rangle} \sigma_{ \Pi_v + tv}(xz_0)}{x^d}
  (-1)^d \sum_{k \geq 0}  
  B_k(t\langle v,z_0\rangle,\va_v) \frac{x^k}{k!}  
     \label{latter general contributions}
\end{align}
valid for almost all $z_0 \in \C^d$.
In the penultimate step we used the definition of the Barnes  polynomials \eqref{Barnes}:
\begin{equation} \label{def. of Barnes yet again}
\frac{ x^d e^{t x} }{      (e^{a_1 x} - 1) \cdots (e^{a_d x} - 1)    } :=
\sum_{k \geq 0}  B_k(t, \va) \frac{x^k}{k!},
\end{equation}
where we replaced $a_k$ by  $ \langle w_k(v), z_0 \rangle$. 
In other words, we used  the definition of the Barnes polynomial, but with the vector 
 $\va_{v}:=  \Big( \langle  w_1(v), z_0 \rangle, \dots, \langle w_d(v), z_0 \rangle \Big)^T$.
To summarize, by using \eqref{first general step} we may sum the latter contributions \eqref{latter general contributions}, one from each of
 the vertex tangent cones, to obtain:
\begin{align} \label{right side in step 1 of theorem 1}
 \sigma_{t\P}(xz_0) 
 &= 
 \sigma_{t\K_{v_1}}(xz_0)  + \cdots +  \sigma_{t\K_{v_N}}(xz_0) \notag\\  
&=  
\sum_{v \in \rm{V}} 
 e^{-t x\langle v, z_0 \rangle} 
 \sigma_{ \Pi_v + tv}(xz_0)
  (-1)^d
 \sum_{k \geq 0}  B_k(t  \langle v, z_0 \rangle,   \va_{v}) \frac{x^{k-d}}{k!},
\end{align}
valid for almost all $z_0 \in \C^d$.

\medskip
\noindent 
{\bf Step $2$}. Next, we massage a  bit the contribution of  
$e^{-t x\langle v, z_0 \rangle} 
 \sigma_{ \Pi_v + tv}(xz_0)$  as well:
\begin{align}
    e^{-t x\langle v, z_0 \rangle} 
 \sigma_{\Pi_v + tv}(xz_0) 
 &=
 e^{x\langle -tv, z_0 \rangle}
 \sum_{p \in (\Pi_v +tv) \cap \Z^d}  
 e^{x \langle p, z_0 \rangle} \notag \\
 &=
 \sum_{p \in (\Pi_v +tv) \cap \Z^d}  
 e^{x \langle p-tv, z_0 \rangle} \notag\\
 &=
 \sum_{q \in (\Z^d - tv) \cap \Pi_v}  
 e^{x \langle q, z_0 \rangle} 
\label{pre lattice flow 1} \\
&=
 \sum_{q \in \Z^d - tv \; (  \hspace{-.1cm}   \bmod \Pi_v)}  
 e^{x \langle q, z_0 \rangle}, 
 \label{lattice flow 1}
\end{align}
where we have made the change of variable $q:= p - tv$ in the penultimate step 
\eqref{pre lattice flow 1} above. This change of variable allows us to think of the integer point transform in \eqref{lattice flow 1} as a {\bf lattice flow on a compact torus}, as defined above in \eqref{def:lattice flow}.  Namely, as the positive real $t$ increases, the points $q$ range over $\Z^d - tv \mod \Pi_v$, inside the fixed torus defined by the parallelepiped $\Pi_v$. 
Inserting \eqref{lattice flow 1} into 
\eqref{right side in step 1 of theorem 1}, we now have:
\begin{align}
\sigma_{t\P}(xz_0) \label{to be used in vanishing identities proof}
&=   (-1)^d \sum_{v \in V} 
  \sum_{q \in  \Z^d - tv \; (  \hspace{-.1cm}   \bmod \Pi_v)}  
 e^{x \langle q, z_0 \rangle}
\sum_{k \geq 0}  B_k(t  \langle v, z_0 \rangle,   \va_{v}) \frac{x^{k-d}}{k!} \\
&= (-1)^d \sum_{v \in V}    
\sum_{q \in  \Z^d - tv \; (  \hspace{-.1cm}   \bmod \Pi_v)}      
\sum_{n\geq 0}   {\langle  q, z_0 \rangle}^n    \frac{x^n}{n!} 
\sum_{k \geq 0}  B_k(t  \langle v, z_0 \rangle,   \va_{v}) \frac{x^{k-d}}{k!} \\
&=  (-1)^d    \sum_{v \in V}                   \label{final sum2, theorem 2}
\sum_{q \in  \Z^d - tv \; (  \hspace{-.1cm}   \bmod \Pi_v)}     
 \sum_{n, k \geq 0}    
 {\langle  q, z_0 \rangle}^n 
  B_k(t  \langle v, z_0 \rangle, \va_{v})
   \frac{x^{k+n-d}}{k! n!}.
\end{align} 

\medskip
\noindent
{\bf Step $3$}. We now compare the constant terms in x, on both sides of  
\eqref{final sum2, theorem 2}.  On its left-hand side, the constant term is $\sigma_{t\P}(0) = L_\P(t)$, by definition.  On the right-hand side of \eqref{final sum2, theorem 2}, we need to collect all terms whose indices satisfy the additional constraint $k+n -d=0$:
\begin{align}
L_\P(t) 
& =   
\frac{(-1)^d}{d!}    
\sum_{v \in V}                       
 \sum_{q \in  \Z^d - tv \; (  \hspace{-.1cm}   \bmod \Pi_v)}    
 \sum_{k=0}^d  
 \binom{d}{k} 
 {\langle  q, z_0 \rangle}^{d-k} 
  B_k(t  \langle v, z_0 \rangle,   \va_{v}) \\  \label{first form of Ehrhart poly 2}
&=  
  \frac{(-1)^d}{d!}    \sum_{v \in V}      
  \sum_{k=0}^d  
  \binom{d}{k}  
  B_k(t  \langle v, z_0 \rangle,   \va_{v})
  \sum_{q \in  \Z^d - tv \; (  \hspace{-.1cm}   \bmod \Pi_v)}   
   {\langle  q, z_0 \rangle}^{d-k},
\end{align} 
valid for almost all $z_0 \in \C^d$, and for every rational $t>0$.

\medskip
\noindent
{\bf Step $4$}.
It remains to extend \eqref{first form of Ehrhart poly 2} from positive rational values of $t$ to all real $t>0$. Fix a generic $z_0\in\C^d$, and denote the right-hand side of \eqref{first form of Ehrhart poly 2} by $\mathcal {\rm Right}_{\P}(t, z_0)$.

Because of the rationality of $\P$, its vertices, the matrices $W_v$, and $t$, both the lattice-point set
$t\P\cap\Z^d$
and, for each vertex $v$, the finite set
$
(tv+\Pi_v)\cap\Z^d
$
can change only at a locally finite set of rational values of $t$. Hence, on every connected component $J$ of the complement of these breakpoints, $L_{\P}(t)$ is constant, while 
${\rm Right}_{\P}(t, z_0)$ is a polynomial function of $t$. Because the equality 
$
L_{\P}(t)={\rm Right}_{\P}(t, z_0)
$
has already been proved for every positive rational $t$,  the two sides agree throughout $J$. At the breakpoints themselves, which are rational, the identity follows directly from the positive-rational case. Thus \eqref{first form of Ehrhart poly 2} holds for every real $t>0$.


\bigskip
To prove part 
\eqref{part b of Theorem 2}, we expand the Barnes polynomial in the latter formula, in terms of the Barnes numbers.  Namely, we use
$B_k(t, \va) = \sum_{j=0}^k  
\binom{k}{j}
 B_{j}(\va) t^{k-j}$, in equation 
 \eqref{first form of Ehrhart poly 2}:
\begin{align} 
\label{the mess before the final step}
L_\P(t) 
& =    \frac{(-1)^d}{d!}    
\sum_{v \in V}                          
    \sum_{k=0}^d   
    \sum_{j=0}^k   
\binom{d}{k} \binom{k}{j}
 B_{j}(\va_v) t^{k-j}  \langle v, z_0 \rangle^{k-j}
   \sum_{q \in  \Z^d - tv \; (  \hspace{-.1cm}   \bmod \Pi_v)}   
   {\langle  q, z_0 \rangle}^{d-k},
\end{align} 
To isolate the coefficient of $t^r$, we retain precisely those terms for which
$k-j=r$.  Thus $k=r+j$, and the constraints $0\leq j\leq k\leq d$
become $0\leq j\leq d-r$.  Therefore \eqref{the mess before the final step}
implies that
\begin{align}
c_r(t)
&=
\frac{(-1)^d}{d!}
\sum_{v\in V}
\sum_{j=0}^{d-r}
\binom{d}{r+j}\binom{r+j}{j}
B_j(\va_v)\langle v,z_0\rangle^r
\sum_{q\in\Z^d-tv\;(\bmod\Pi_v)}
\langle q,z_0\rangle^{d-r-j}
\notag\\
&=
\frac{(-1)^d}{d!}
\sum_{v\in V}
\sum_{j=0}^{d-r}
\binom{d}{r,j}
B_j(\va_v)\langle v,z_0\rangle^r
\sum_{q\in\Z^d-tv\;(\bmod\Pi_v)}
\langle q,z_0\rangle^{d-r-j},
\end{align}
valid for almost all $z_0\in\C^d$.  This proves part
\eqref{part b of Theorem 2}.
\end{proof}


\section{Proof of Theorem  
\ref{thm:Discrete moments of rational polytopes}}

We recall the simple observation that the integer point transform of $\P$ is a generating function for the discrete moments over $\P$. Precisely, we have:
\begin{equation}
\label{starting point for discrete moments proof}
    \sigma_{t\P}(xz_0):=
\sum_{n \in t\P \cap \Z^d}    
e^{ \langle n,  xz_0 \rangle}
=
\sum_{n \in t\P \cap \Z^d}    
\sum_{k\geq 0} \frac{x^k}{k!} 
{\langle n,  z_0 \rangle}^k
:=
\sum_{k\geq 0} \frac{x^k}{k!}
\mu_k(t\P, z_0).
\end{equation}

\begin{proof}
\hypertarget{Discrete moments of polytopes}
{(of Theorem \ref{thm:Discrete moments of rational polytopes})}
We begin with equation 
\eqref{final sum2, theorem 2} from the proof of Theorem \ref{thm:Ehrhart quasi-polynomial from Barnes}:

\begin{align}
\label{begin with this, from proof of theorem 2}
\sigma_{t\P}(xz_0)
&=  (-1)^d    \sum_{v \in V}    
\sum_{q \in  \Z^d - tv \; (  \hspace{-.1cm}   \bmod \Pi_v)} 
\sum_{n, k \geq 0}    
 {\langle  q, z_0 \rangle}^n 
  B_k(t  \langle v, z_0 \rangle,   \va_{v})
   \frac{x^{k+n-d}}{k! n!}.
\end{align} 
We fix any integer $m>0$, and we compare the coefficients of $x^m$ on
both sides of  \eqref{begin with this, from proof of theorem 2}.
Using \eqref{starting point for discrete moments proof}, we trivially see that the coefficient of $x^m$ on the left-hand side of
\eqref{begin with this, from proof of theorem 2} is
$\frac{1}{m!}\mu_m(t\P, z_0)$.

On the right-hand side of \eqref{begin with this, from proof of theorem 2}, we need to collect all terms whose indices satisfy the  constraint 
$k+n -d=m$:
\begin{align}
  \frac{1}{m!}\mu_m(t\P, z_0) 
&= 
(-1)^d    
\sum_{v \in V}                  
\sum_{q \in  \Z^d - tv \; (  \hspace{-.1cm}   \bmod \Pi_v)} 
\sum_{k=0}^{d+m}  
\frac{1}{k! (d+m-k)!} 
 {\langle  q, z_0 \rangle}^{d+m-k} 
  B_k(t  \langle v, z_0 \rangle,   \va_{v}), \label{set-up for part b}
\end{align} 
valid for almost all $z_0 \in \C^d$. This completes the proof of part 
\eqref{part a of discrete moments}.

\noindent
{\bf Step $3$}. 
To prove part 
\eqref{part b of Theorem 2}, 
we expand the Barnes polynomial in the latter formula, in terms of the Barnes numbers.  Namely, we use
$B_k(t, \va) = \sum_{j=0}^k  
\binom{k}{j}
 B_{j}(\va) t^{k-j}$, in equation \eqref{set-up for part b}:
 \eqref{first form of Ehrhart poly 2}:

\begin{align} 
\label{the new mess before the final step}
&  \tfrac{(-1)^d }{m!}\mu_m(t\P, z_0)  
 =      
\sum_{v \in V}                          
    \sum_{k=0}^{d +m} 
    \sum_{j=0}^k   
     \tfrac{1}{k! (d+m-k)!} 
 \binom{k}{j}
 B_{j}(\va_v) t^{k-j}  \langle v, z_0 \rangle^{k-j}
   \sum_{q \in  \Z^d - tv \; (  \hspace{-.1cm}   \bmod \Pi_v)}   
   {\langle  q, z_0 \rangle}^{d+m-k},
\end{align} 
Finally, we just have to sum over all relevant $k, j$, such that  $k-j := r$, with $r$ fixed, 
so that we may isolate 
the coefficient of $t^r$ on the right-hand side. 
 We use the following elementary combinatorial identity, when setting $k-j = r$, which holds for any function $F(j,k)$ on $\Z_{\geq 0} \times \Z_{\geq 0}$:
\[
\sum_{k=0}^{d+m} \sum_{j= 0}^k F(j,k) = \sum_{j=0}^{d+m-r} F(j, r + j).
\]
In our context,
\eqref{the new mess before the final step} therefore implies that the coefficient
$d_r(t)$ of $t^r$ is
\begin{align}
d_r(t)
&=
(-1)^d m!
\sum_{v\in V}
\sum_{j=0}^{d+m-r}
\frac{B_j(\va_v)\langle v,z_0\rangle^r}
{r!\,j!\,(d+m-r-j)!}
\sum_{q\in\Z^d-tv\;(\bmod\Pi_v)}
\langle q,z_0\rangle^{d+m-r-j}
\notag\\
&=
\frac{(-1)^d m!}{(d+m)!}
\sum_{v\in V}
\sum_{j=0}^{d+m-r}
\binom{d+m}{r,j}
B_j(\va_v)\langle v,z_0\rangle^r
\sum_{q\in\Z^d-tv\;(\bmod\Pi_v)}
\langle q,z_0\rangle^{d+m-r-j},
\label{eq:corrected-moment-quasi-coefficient-proof}
\end{align}
valid for almost all $z_0\in\C^d$.

To prove part~\eqref{Extension of Linke's ODE}, fix
$t>0$ away from the locally finite set of breakpoints at which one of
the points in $\bigl(\Pi_v+tv\bigr)\cap\Z^d$ changes.  In a neighborhood of such a
$t$, every one of these finite sets is constant.  Rewriting the final
sum in \eqref{eq:corrected-moment-quasi-coefficient-proof} as
\[
\sum_{q\in(\Pi_v+tv)\cap\Z^d}
\bigl(\langle q,z_0\rangle-t\langle v,z_0\rangle\bigr)^{d+m-r-j},
\]
we may differentiate term by term.  The term with
$j=d+m-r$ has exponent zero and contributes nothing, so
\begin{align*}
\frac{d}{dt}d_r(t)
&=
-(-1)^d m!
\sum_{v\in V}
\langle v,z_0\rangle^{r+1}
\sum_{j=0}^{d+m-r-1}
\frac{B_j(\va_v)}
{r!\,j!\,(d+m-r-j-1)!}\\
&\hspace{2.5cm}\times
\sum_{q\in(\Pi_v+tv)\cap\Z^d}
\bigl(\langle q,z_0\rangle-t\langle v,z_0\rangle\bigr)^{d+m-r-j-1}\\
&=-(r+1)d_{r+1}(t),
\end{align*}
where the final equality follows by comparing the preceding expression
with \eqref{eq:corrected-moment-quasi-coefficient-proof} with $r$
replaced by $r+1$.  This proves the differential equation for every
$t>0$ away from the locally finite set of breakpoints.
\end{proof}


\section{Proof of Theorem \ref{thm:MainVanishingIdentities}:
vanishing identities for rational simple polytopes}

\begin{proof} 
\hypertarget{Proof of Theorem 3}
{(of Theorem \ref{thm:MainVanishingIdentities})}
Using equation 
\eqref{final sum2, theorem 2}
from the proof of Theorem \ref{thm:Ehrhart.from.Bernoulli.Barnes}, 
we recall:
\begin{align} \label{just before considering singularities}
\sigma_{t\P}(x z_0) 
&= (-1)^d    \sum_{v \in V}            
\sum_{q \in  \Z^d - tv \; (  \hspace{-.1cm}   \bmod \Pi_v)}     
\sum_{n, k \geq 0}  {\langle  q, z_0 \rangle}^n 
  B_k(t  \langle v, z_0 \rangle, \va_{v})
   \frac{x^{k+n-d}}{k! n!}
\end{align}
To prove part \eqref{part a of vanishing identities}, we begin by noting that the left-hand-side of the latter identity is an entire function of $x$.   We therefore  know that for each fixed $n<0$, the coefficient
of any $x^{n}$ on the right-hand-side must be identically zero.  
We now compare the coefficient of $x^{m-d}$ on both sides, and for each fixed value of $m \in \{0, \dots, d-1\}$, we must get vanishing identities. Precisely, we get the following $d$ identities for each rational simple polytope in $\R^d$, by summing 
 over all $k, n \geq 0$ in \eqref{just before considering singularities} 
 with $k+n= m$ fixed:
\begin{align}\label{identities1}
0 
&=   \sum_{k=0}^m 
\binom{m}{k}   
\sum_{v \in V}      
  B_k \Big(t  \langle v, z_0 \rangle,   \va_{v} \Big)   
\sum_{q \in  \Z^d - tv \; (  \hspace{-.1cm}   \bmod \Pi_v)}       
{\langle  q, z_0 \rangle}^{m-k}.
\end{align}  

\noindent
To prove part \eqref{part b of vanishing identities}, we expand the Barnes polynomial as
\[
B_k\!\left(t\langle v,z_0\rangle,\va_v\right)
=
\sum_{j=0}^k
\binom{k}{j}
B_j(\va_v)
\langle v,z_0\rangle^{k-j}t^{k-j}.
\]
Substituting this expansion into part~(a) and collecting terms with the
same power \(t^r\), we obtain
\[
\sum_{r=0}^m t^r C_r(t)=0,
\]
where
\[
C_r(t)
=
\frac{m!}{r!}
\sum_{v\in V}
\langle v,z_0\rangle^r
\sum_{j=0}^{m-r}
\frac{B_j(\va_v)}{j!(m-r-j)!}
\sum_{q\in\Z^d-tv\;(\bmod\Pi_v)}
\langle q,z_0\rangle^{m-r-j}.
\]
The functions \(C_r(t)\) are not constant in \(t\), so we cannot yet
compare coefficients directly. However, since all vertices of
\(\mathcal P\) are rational, there exists a positive integer \(T\) such that
\[
Tv\in\Z^d
\qquad
\text{for every }v\in V.
\]
It follows that every lattice-flow moment occurring in \(C_r(t)\) is
\(T\)-periodic: $C_r(t+T)=C_r(t)$.

We now fix \(t_0>0\), so that the periodic functions \(C_r(t)\) are
frozen at their values \(C_r(t_0)\), turning the identity into an ordinary
polynomial identity in \(t_0+\ell T\) whose coefficients can be compared.
Indeed, evaluating the preceding identity at
$t=t_0+\ell T$, with $\ell\in\Z_{\geq 0}$, 
and using periodicity, we obtain
\[
\sum_{r=0}^m
(t_0+\ell T)^r C_r(t_0)=0
\qquad
\text{for every }\ell\in\Z_{\geq 0}.
\]
Thus the polynomial $P_{t_0}(x) := \sum_{r=0}^m C_r(t_0)x^r$
vanishes at each of the infinitely many distinct points $x=t_0+\ell T$, with  $\ell\in\Z_{\geq 0}$.
We conclude that  \(P_{t_0}\) is the zero polynomial, and consequently
$C_r(t_0)=0$
for every $0\leq r\leq m$, which 
is precisely part~(b).
\end{proof}


\section{Proof of Theorem \ref{thm:moment-reciprocity}}

\begin{proof}
\hypertarget{Proof of Theorem: moment-reciprocity}
{(of Theorem \ref{thm:moment-reciprocity})}
Recalling the definition of \(\mathbf{a}_v\) from \eqref{def:a_v},
we write
\[
\mathbf{a}_v
:=
(a_{v,1},\ldots,a_{v,d})^T,
\]
where $a_{v,j}
:=
\langle w_j(v),z\rangle$.
We fix \(m\in\mathbb{Z}_{\geq 0}\), \(t>0\), and a generic vector
\(z\in\mathbb{C}^d\). 
For each vertex \(v\in V\), we define
$
A_v 
:= 
a_{v,1}+\cdots+a_{v,d} 
=
\left\langle
w_1(v)+\cdots+w_d(v),z
\right\rangle$, 
and 
\begin{equation}
\Pi_v^{*}
:=
\left\{
\sum_{j=1}^{d}\lambda_jw_j(v)
\,:\,
0<\lambda_j\leq 1
\right\}
=
w_1(v)+\cdots+w_d(v)-\Pi_v.
\end{equation}
We let
$
K_{tv}^{o}
:=
tv+\sum_{j=1}^{d}\mathbb{R}_{>0}w_j(v)
$
be the open tangent cone at the vertex \(tv\). The open-cone analogue of the
local rational-function formula is
\begin{equation}
\begin{aligned}
K_{tv}^{o}\cap\mathbb{Z}^d
&=
\bigsqcup_{p\in(tv+\Pi_v^{*})\cap\mathbb{Z}^d}
\left(
p+\sum_{j=1}^{d}\mathbb{Z}_{\geq 0}w_j(v)
\right),\\
\sigma_{K_{tv}^{o}}(xz)
&=
\frac{
\sigma_{\,tv+\Pi_v^{*}}(xz)
}{
\displaystyle\prod_{j=1}^{d}
\left(1-e^{xa_{v,j}}\right)
}.
\end{aligned}
\label{eq:moment-reciprocity-open-cone-transform}
\end{equation}
Indeed, every \(\alpha_j>0\) has a unique decomposition
\(\alpha_j=n_j+\lambda_j\), where
\(n_j\in\mathbb{Z}_{\geq 0}\) and \(0<\lambda_j\leq 1\).
This gives the disjoint decomposition in the first line of
\eqref{eq:moment-reciprocity-open-cone-transform}, and hence the
geometric-series identity in the second line.
Now we'll use the open form of Brion's identity:
\begin{equation}
\sigma_{t\mathcal{P}^{o}}(xz)
=
\sum_{v\in V}\sigma_{K_{tv}^{o}}(xz).
\label{eq:moment-reciprocity-open-brion}
\end{equation}
The exponential integer-point transform of the open polytope encodes
its discrete moments:
\begin{equation}
\begin{aligned}
\sigma_{t\mathcal{P}^{o}}(xz)
&=
\sum_{p\in t\mathcal{P}^{o}\cap\mathbb{Z}^d}
e^{x\langle p,z\rangle}\\
&=
\sum_{p\in t\mathcal{P}^{o}\cap\mathbb{Z}^d}
\sum_{r\geq 0}
\langle p,z\rangle^r\frac{x^r}{r!}\\
&=
\sum_{r\geq 0}
M_r^{o}(\mathcal{P};t,z)\frac{x^r}{r!},
\qquad
M_m^{o}(\mathcal{P};t,z)
=
m![x^m]\sigma_{t\mathcal{P}^{o}}(xz).
\end{aligned}
\label{eq:open-moment-generating-function}
\end{equation}

We next expand each local open-cone integer transform. From
\eqref{eq:moment-reciprocity-open-cone-transform}, we have
\begin{equation}
\begin{aligned}
\sigma_{K_{tv}^{o}}(xz)
&=
\frac{
\sigma_{\,tv+\Pi_v^{*}}(xz)
}{
\displaystyle\prod_{j=1}^{d}
\left(1-e^{xa_{v,j}}\right)
}\\
&=
\frac{
e^{-x\left(t\langle v,z\rangle+A_v\right)}
\sigma_{\,tv+\Pi_v^{*}}(xz)
}{x^d}
\,
\frac{
x^d e^{x\left(t\langle v,z\rangle+A_v\right)}
}{
\displaystyle\prod_{j=1}^{d}
\left(1-e^{xa_{v,j}}\right)
}.
\end{aligned}
\label{eq:moment-reciprocity-open-factorization}
\end{equation}
The finite numerator in the first factor satisfies
\begin{align}
e^{-x\left(t\langle v,z\rangle+A_v\right)}
\sigma_{\,tv+\Pi_v^{*}}(xz)
&=
\sum_{p\in(tv+\Pi_v^{*})\cap\mathbb{Z}^d}
e^{
x\left(
\langle p,z\rangle
-t\langle v,z\rangle
-A_v
\right)
}\\
&\qquad=
\sum_{p\in(tv+\Pi_v^{*})\cap\mathbb{Z}^d}
e^{-x\left\langle
tv+w_1(v)+\cdots+w_d(v)-p,z
\right\rangle}\\
&\qquad=
\sum_{q\in(\mathbb{Z}^d+tv)\cap\Pi_v}
e^{-x\langle q,z\rangle}.
\label{eq:moment-reciprocity-open-numerator}
\end{align}
The final equality follows from the bijection
\[
p
\longmapsto
q
:=
tv+w_1(v)+\cdots+w_d(v)-p.
\]
So we see that since $\Pi_v^{*}=
w_1(v)+\cdots+w_d(v)-\Pi_v$,
the condition 
\(p-tv\in\Pi_v^{*}\) is equivalent to \(q\in\Pi_v\).
Moreover, because
$p,\, w_1(v)+\cdots+w_d(v) \in\mathbb{Z}^d$, 
we have \(q\in\mathbb{Z}^d+tv\).
For the second factor in
\eqref{eq:moment-reciprocity-open-factorization}, we use the defining
generating function of the Barnes polynomials, followed by the Barnes
reciprocity formula
\begin{equation}
(-1)^kB_k(-u,\mathbf{a})
=
B_k(u+a_1+\cdots+a_d,\mathbf{a}),
\end{equation}
given in Appendix \ref{Basic properties of Barnes}, property
\eqref{reciprocity}, with
$A_v:= a_1+\cdots+a_d$.
It follows that
\begin{align}
\label{eq:moment-reciprocity-barnes-factor}
\frac{
x^d e^{x\left(t\langle v,z\rangle+A_v\right)}
}{
\displaystyle\prod_{j=1}^{d}
\left(1-e^{xa_{v,j}}\right)
}
&=
(-1)^d
\frac{
x^d e^{x\left(t\langle v,z\rangle+A_v\right)}
}{
\displaystyle\prod_{j=1}^{d}
\left(e^{xa_{v,j}}-1\right)
}\\
&\qquad=
(-1)^d
\sum_{k\geq 0}
B_k\!\left(
t\langle v,z\rangle+A_v,\mathbf{a}_v
\right)
\frac{x^k}{k!}\\
&\qquad=
(-1)^d
\sum_{k\geq 0}
(-1)^k
B_k\!\left(
-t\langle v,z\rangle,\mathbf{a}_v
\right)
\frac{x^k}{k!}.
\end{align}
Combining
\eqref{eq:moment-reciprocity-open-factorization},
\eqref{eq:moment-reciprocity-open-numerator}, and
\eqref{eq:moment-reciprocity-barnes-factor}, and then expanding the
exponential, gives
\begin{equation}
\begin{aligned}
\sigma_{K_{tv}^{o}}(xz)
&=
(-1)^d
\sum_{q\in(\mathbb{Z}^d+tv)\cap\Pi_v}
e^{-x\langle q,z\rangle}
\sum_{k\geq 0}
(-1)^k
B_k\!\left(
-t\langle v,z\rangle,\mathbf{a}_v
\right)
\frac{x^{k-d}}{k!}\\
&=
(-1)^d
\sum_{q\in(\mathbb{Z}^d+tv)\cap\Pi_v}
\left(
\sum_{n\geq 0}
(-1)^n
\langle q,z\rangle^n
\frac{x^n}{n!}
\right)\\
&\qquad\qquad\times
\left(
\sum_{k\geq 0}
(-1)^k
B_k\!\left(
-t\langle v,z\rangle,\mathbf{a}_v
\right)
\frac{x^{k-d}}{k!}
\right)\\
&=
(-1)^d
\sum_{q\in(\mathbb{Z}^d+tv)\cap\Pi_v}
\sum_{n,k\geq 0}
(-1)^{n+k}
\langle q,z\rangle^n
B_k\!\left(
-t\langle v,z\rangle,\mathbf{a}_v
\right)
\frac{x^{n+k-d}}{n!\,k!}.
\end{aligned}
\label{eq:moment-reciprocity-local-laurent-expansion}
\end{equation}

We now sum over the vertices using
\eqref{eq:moment-reciprocity-open-brion} and extract the coefficient
of \(x^m\). The relevant condition is
$n+k-d=m$.
We therefore have:
\begin{equation}
\begin{aligned}
M_m^{o}(\mathcal{P};t,z)
&=
m![x^m]\sigma_{t\mathcal{P}^{o}}(xz)\\
&=
(-1)^d m!
\sum_{v\in V}
\sum_{q\in(\mathbb{Z}^d+tv)\cap\Pi_v}
\sum_{\substack{n,k\geq 0\\ n+k=d+m}}
(-1)^{n+k}
\frac{
\langle q,z\rangle^n
B_k\!\left(
-t\langle v,z\rangle,\mathbf{a}_v
\right)
}{
n!\,k!
}\\
&=
(-1)^m m!
\sum_{v\in V}
\sum_{k=0}^{d+m}
\frac{
B_k\!\left(
-t\langle v,z\rangle,\mathbf{a}_v
\right)
}{
k!\,(d+m-k)!
}
\sum_{q\in(\mathbb{Z}^d+tv)\cap\Pi_v}
\langle q,z\rangle^{d+m-k}.
\end{aligned}
\label{eq:interior-moment-barnes-formula}
\end{equation}

On the other hand, substituting 
\(s=-t\) into
\eqref{eq:barnes-extension-of-moments} gives
\begin{align}
(-1)^{d+m}M_m(\mathcal{P};-t,z)
&=
(-1)^m m!
\sum_{v\in V}
\sum_{k=0}^{d+m}
\frac{
B_k\!\left(
-t\langle v,z\rangle,\mathbf{a}_v
\right)
}{
k!\,(d+m-k)!
}
\sum_{q\in(\mathbb{Z}^d+tv)\cap\Pi_v}
\langle q,z\rangle^{d+m-k}\\
&=
M_m^{o}(\mathcal{P};t,z).
\end{align}
The final equality is precisely
\eqref{eq:interior-moment-barnes-formula}, proving
\eqref{eq:moment-reciprocity-statement} for generic \(z\).

Finally, the left-hand side of
\eqref{eq:moment-reciprocity-statement} is a homogeneous polynomial
of degree \(m\) in \(z\). Since the identity holds outside the finite
union of hyperplanes on which some
\(\langle w_j(v),z\rangle\) vanishes, the total Barnes expression on
the right has a unique polynomial extension across those
hyperplanes. The reciprocity identity therefore extends to every
\(z\in\mathbb{C}^d\).
\end{proof}


\section{Proof of Corollary \ref{cor:general rational polytopes}}

\begin{proof} 
\hypertarget{Proof of Corollary:extension of main theorem to rational polytopes}
{(of Corollary \ref{cor:general rational polytopes})}
For each vertex \(v\in V\), we let
$C_v:=\K_v-v$.

{\bf Step} $1$.  We triangulate each vertex tangent cone $C_v$ into partially-open simplicial cones.
To do so rigorously, we choose a rational simplicial triangulation
\[
C_v=\bigcup_{j=1}^{N_v}\overline C_{v,j}.
\]
Throughout this proof, a maximal cone of a triangulation means a  \(d\)-dimensional simplicial cone. We write
\(\overline{C}_{v,j}\) for the closed maximal cones of the triangulation
and \(C_{v,j}\) for their partially open versions obtained by the generic
perturbation. Thus, the closed cones may intersect along proper faces,
whereas the partially open cones form a disjoint decomposition of \(C_v\).

 A standard generic
perturbation assigns every common face to exactly one maximal cone, as follows.  We choose
\(\xi_v\in\operatorname{int}(C_v)\) outside the linear span of every facet
appearing in the triangulation, and define
\begin{equation}
\label{eq:half-open-triangulation-general-polytope}
C_{v,j}
:=
\left\{
y\in\overline C_{v,j}:
y+\varepsilon\xi_v\in\overline C_{v,j}
\text{ for all sufficiently small }\varepsilon>0
\right\}.
\end{equation}
It follows that
$C_v = \bigsqcup_{j=1}^{N_v}C_{v,j}$,
and that the \(C_{v,j}\) are partially-open simplicial cones that form 
a disjoint decomposition of \(C_v\).

{\bf Step} $2$. \ Let \(w_1(v,j),\ldots,w_d(v,j)\) be the primitive integer edges of \(C_{v,j}\), and let \(\Pi_{v,j}\) be its half-open
 parallelepiped, with the boundary convention prescribed in
Lemma~\ref{lem:extension to any rational polytope}. That lemma gives
\begin{equation}
\label{eq:partially-open-cone-transform-general-polytope}
\sigma_{tv+C_{v,j}}(xz)
=
\frac{\sigma_{tv+\Pi_{v,j}}(xz)}
{\displaystyle
\prod_{i=1}^d
\left(
1-e^{x\langle w_i(v,j),z\rangle}
\right)},
\end{equation}
as an identity of rational functions.
We fix \(t\in\Q_{>0}\), and choose \(z\in\C^d\) such that
\(\langle w_i(v,j),z\rangle\neq0\) for every \(v,j,i\). By Brion's
theorem and the disjoint decomposition $C_v = \bigsqcup_{j=1}^{N_v}C_{v,j}$, we have:
\begin{equation}
\label{eq:Brion-disjoint-general-polytope}
\sigma_{t\P}(xz)
=
\sum_{v\in V}
\sum_{j=1}^{N_v}
\sigma_{tv+C_{v,j}}(xz).
\end{equation}
Moreover,
\begin{equation}
\label{eq:recentered-numerator-general-polytope}
e^{-tx\langle v,z\rangle}
\sigma_{tv+\Pi_{v,j}}(xz)
=
\sum_{p\in\Z^d-tv\;(\bmod\Pi_{v,j})}
e^{x\langle p,z\rangle}.
\end{equation}

{\bf Step} $3$.  \ Applying the Barnes generating function to the denominator in
\eqref{eq:partially-open-cone-transform-general-polytope}, expanding
\eqref{eq:recentered-numerator-general-polytope} in powers of \(x\),
and taking the constant Laurent coefficient gives
\[
L_\P(t)
=
\frac{(-1)^d}{d!}
\sum_{v\in V}
\sum_{j=1}^{N_v}
\sum_{k=0}^d
\binom{d}{k}
B_k\left(
t\langle v,z\rangle,\va_{v,j}
\right)
\sum_{p\in\Z^d-tv\;(\bmod\Pi_{v,j})}
\langle p,z\rangle^{d-k}.
\]
This is exactly the calculation in the proof of
Theorem~\ref{thm:Ehrhart quasi-polynomial from Barnes}, but now applied to
each partially-open simplicial cone \(C_{v,j}\). So this proves the latter formula
for every \(t\in\Q_{>0}\).

{\bf Step} $4$.  
 For rigorous completeness, the extension to real \(t>0\) follows from the same
chamber argument used in Theorem 
\ref{thm:Ehrhart quasi-polynomial from Barnes}. 
The left-hand side can change
only when a lattice point meets a moving rational facet of \(t\P\).
 All isolated breakpoints are rational, and they
form a locally finite set. Between consecutive breakpoints the
left-hand side is constant and the right-hand side is polynomial in
\(t\); since they agree at every rational \(t\), they agree throughout
each such interval. Equality at the breakpoints follows from the
already established rational case. Hence the formula holds for every
real \(t>0\).

Finally, since \(W_{v,j}^{-1}v\in\Q^d\) (see Observation \ref{Observation: explicit lattice flow description} above), the corresponding lattice
flows have a common positive integral period. Expanding the Barnes
polynomials in their first argument therefore shows that the displayed
expression is a quasi-polynomial.
\end{proof}


\section{Proof of Corollary \ref{Ehrhart for smooth polytopes}: Ehrhart for smooth polytopes}

\begin{proof} 
\hypertarget{Proof of smooth formula}
{(of Corollary \ref{Ehrhart for smooth polytopes})}
The smooth polytope assumption means that for all vertices $v\in \P$, there is exactly one  integer point in the vertex parallelepiped $\Pi_v$, namely the origin. Therefore the discrete moments, for each $0\leq k \leq d-1$, degenerate to: 
\begin{equation}
    \sum_{p \in \Pi_v \cap \Z^d}   
{\langle  p, z \rangle}^{d-k} 
= 0,
\end{equation}
for all vertices $v \in \P$.  
Moreover, for $k=d$, we have 
\[
\sum_{p \in \Pi_v \cap \Z^d}   
{\langle  p, z \rangle}^{d-k}
=\sum_{p \in \Pi_v \cap \Z^d} 1 = 1,
\] 
so that we have:
\begin{align*}
  L_\P(t)  
&=     
\frac{(-1)^d}{d!} 
\sum_{k=0}^d  \binom{d}{k}
\sum_{v \in V}                             
B_k(t  \langle v, z \rangle,   \va_{v})
   \sum_{p \in \Pi_v \cap \Z^d}   
   {\langle  p, z \rangle}^{d-k}\\
&=
\frac{(-1)^d}{d!} 
\sum_{v \in V}                             
B_d(t  \langle v, z \rangle,   \va_{v}).
\end{align*}
\end{proof}


\section{More examples, and some identities for special cases of some of the Theorems}
\label{sec: examples}

 \bigskip 
 \begin{example}  \label{first example}
 In dimension $d=2$,
 Corollary \ref{constant term 1 identity for smooth polytopes}
 tells us that for any smooth integer polygon $P\subset \R^2$ we have: 
\begin{equation}
1 = \frac{  1}{2}     
 \sum_{v \in \rm{V}}          B_{2}(\va_{v}) 
 =  \frac{  1}{2}    
 \sum_{v \in \rm{V}}    
  \left(      
\frac{1}{6}\left(    \frac{\langle w_2(v), z \rangle}{\langle w_1(v), z \rangle} +
\frac{\langle w_1(v), z \rangle}{\langle w_2(v), z \rangle}   
\right) + \frac{1}{2}
\right),
\end{equation}
for almost all $z \in \C^2$. 
In other words, we have the identity
\begin{equation}
12 =   
 \sum_{v \in \rm{V}}   
 \left(    \frac{\langle w_2(v), z \rangle}{\langle w_1(v), z \rangle} +
\frac{\langle w_1(v), z \rangle}{\langle w_2(v), z \rangle}   +  3
\right),
\end{equation}
valid for almost all $z \in \C^2$.
\hfill $\square$
\end{example}
\begin{example}
\label{Example: dimension 3, the 24-identity}
In dimension $d=3$, 
Corollary \ref{constant term 1 identity for smooth polytopes}
tells us that for any smooth integer polytope $P\subset \R^3$, we have:
\begin{align*} 
1 = \frac{  (-1)^3 }{3!}     
 \sum_{v \in \rm{V}} B_{3}(\va_{v})
= -\frac{1}{6} \sum_{v \in \rm{V}}        
\frac{1}{8a_1(v) a_2(v) a_3(v)}(s_{1}s_{2} - s_{1}^{3}),
\end{align*}
by using \eqref{B_3(a), dim d} of Appendix \ref{Basic properties of Barnes}.
Using the definitions of the power sums $s_j$ and the definition  $a_j(v):= \langle w_j(v), z \rangle$, we arrive at the identity 
\begin{align} 
24 &=  \frac{1}{2}\sum_{v \in \rm{V}}        
\frac{s_{1}^{3}-s_{1}s_{2}  }{a_1(v) a_2(v) a_3(v)}\\
&=
 \frac{1}{2}\sum_{v \in \rm{V}}        
\tfrac{\left(  \langle w_1(v), z \rangle + \langle w_2(v), z \rangle  + \langle w_3(v), z \rangle  \right)^{3}
   - \left(  \langle w_1(v), z \rangle + \langle w_2(v), z \rangle  + \langle w_3(v), z \rangle  \right)
   \left(  \langle w_1(v), z \rangle^2 + \langle w_2(v), z \rangle^2  + \langle w_3(v), z \rangle^2  \right)  }
{\langle w_1(v), z \rangle \langle w_2(v), z \rangle  \langle w_3(v), z \rangle }\\
&=
\sum_{v \in \rm{V}}
\left(
\tfrac{ \langle w_1(v), z\rangle}{\langle w_2(v), z\rangle}
+\tfrac{ \langle w_2(v), z\rangle}{\langle w_1(v), z\rangle}
+\tfrac{ \langle w_2(v), z\rangle}{\langle w_3(v), z\rangle}
+\tfrac{ \langle w_3(v), z\rangle}{\langle w_2(v), z\rangle}
+\tfrac{ \langle w_3(v), z\rangle}{\langle w_1(v), z\rangle}
+\tfrac{ \langle w_1(v), z\rangle}{\langle w_3(v), z\rangle} +3
\right),
\end{align}
valid for almost all $z \in \C^3$.
\hfill $\square$
\end{example}

\subsection{Special cases of Theorem \ref{thm:MainVanishingIdentities}}
We recall the vanishing identities from Theorem \ref{thm:MainVanishingIdentities}:
 \begin{align}
0 &=  \sum_{v \in \rm{V}}    
   \langle v, z \rangle^{r}       
 \sum_{j=0}^{m-r} 
\frac{B_{j}(\va_{v})}{j!(m-r - j)!}    
\sum_{q \in  \Z^d - tv \; (  \hspace{-.1cm}   \bmod \Pi_v)}       
{\langle  q, z \rangle}^{m-r-j},
\end{align} 
for each $(m, r)$ with  $0\leq m \leq  d-1,   \text{  and }  0\leq r \leq m$.

Here we look at some special cases, 
noting that the particular cases of $r = m$ recover known vanishing identities,  by Brion and Vergne \cite{BrionVergne}, while  the other identities appear to be new. 
 Precisely, for $r = m$ with $0 \leq m \leq d-1$, these $d$
 known identities are:
\begin{align}  
\label{Brion vanishing identities}
0 
&=  
\sum_{v \in \rm{V}}    
\langle v, z \rangle^{m}       
B_{0}(\va_{v})           
\sum_{q \in  \Z^d - tv \; (  \hspace{-.1cm}   \bmod \Pi_v)}     1 =
\sum_{v \in V}      \frac{    \langle v, z \rangle^{m}    \vol \Pi_{v}   }{ \prod_{k=1}^d  \langle w_k(v), z \rangle}, 
\end{align}  
 valid for almost all $z\in \C^d$.    Here we used $B_0(\va) = \frac{1}{a_1 \cdots a_d}$, and 
\begin{equation}
\label{basic facts for Pi_v}
\sum_{q \in  \Z^d - tv \; (  \hspace{-.1cm}   \bmod \Pi_v)}  1
=
\displaystyle\sum_{q \in  \Z^d  \; (  \hspace{-.1cm}   \bmod \Pi_v)} 1
= 
|\Pi_v \cap \Z^d|
= 
\vol \Pi_{v},
\end{equation}
for all $t>0$.
The latter equality is a standard fact from the geometry of numbers \cite[Theorem 5.11]{Sinai2}, while the first equality in 
\eqref{basic facts for Pi_v} 
follows from the easy fact that the lattice flow preserves the number of integer points in the compact torus.

 Throughout this section, we fix the dimension $d \geq 2$.
 Regarding the new identities, we first discuss the case $m=1, r=0$, as it presents a special case that is
 of independent interest: equation \eqref{independent interest, m=1, r=0} below. 
 
 \bigskip \noindent
 {\bf Case $(m, r) = (1, 0)$}.
   \  Here Theorem \ref{thm:MainVanishingIdentities} translates into:
\begin{align}   \label{case of m=1, r=0}
0 &=  \sum_{v \in \rm{V}}        
 \sum_{j=0}^{1} 
\frac{B_{j}(\va_{v})}{j!}            
\sum_{q \in  \Z^d - tv \; (  \hspace{-.1cm}   \bmod \Pi_v)}   
\frac{ {\langle  p, z \rangle}^{1 - j} }{(1 - j)!} \\
&=  \sum_{v \in \rm{V}}       
B_{0}(\va_{v})           
\sum_{q \in  \Z^d - tv \; (  \hspace{-.1cm}   \bmod \Pi_v)}   
\langle  p, z \rangle  +      
 \sum_{v \in \rm{V}}           
 B_{1}(\va_{v})           
\sum_{q \in  \Z^d - tv \; (  \hspace{-.1cm}   \bmod \Pi_v)}    1           \\
&=     \sum_{v \in \rm{V}}       \frac{ 1 }{ \prod_{k=1}^d  \langle w_k(v), z \rangle}       
\sum_{q \in  \Z^d - tv \; (  \hspace{-.1cm}   \bmod \Pi_v)}     \langle  p, z \rangle  +      
\sum_{v \in \rm{V}}           B_{1}(\va_{v})   \vol \Pi_{v}   \\  \label{same as Case (1, 0)}
&= \sum_{v \in \rm{V}}       
\frac{ 1 }{ \prod_{k=1}^d  
\langle w_k(v), z \rangle}         
\sum_{q \in  \Z^d - tv \; (  \hspace{-.1cm}   \bmod \Pi_v)}   
\langle  p, z \rangle        
 - \frac{1}{2} \sum_{v \in \rm{V}}     
 \vol \Pi_{v}              
\frac{   \sum_{k=1}^d \langle w_k(v), z \rangle }{\prod_{k=1}^d \langle w_k(v), z \rangle},
\end{align}
so that we have 
\begin{align} \label{independent interest, m=1, r=0}
 \sum_{v \in \rm{V}}       
 \frac{ 1 }{ \prod_{k=1}^d  \langle w_k(v), z \rangle}         
\sum_{q \in  \Z^d - tv \; (  \hspace{-.1cm}   \bmod \Pi_v)}      
\langle  p, z \rangle        
 = \frac{1}{2} \sum_{v \in \rm{V}}     
 \vol \Pi_{v}              
\frac{   \sum_{k=1}^d \langle w_k(v), z \rangle }{\prod_{k=1}^d \langle w_k(v), z \rangle}.
\end{align}


 {\bf Case} $(m, r) = (m, m-1)$, with $1\leq m  \leq d-1$.   \  
Here Theorem \ref{thm:MainVanishingIdentities} gives us  
 \begin{align} 
0 &=  
 \sum_{v \in \rm{V}}    
   \langle v, z \rangle^{r}       
 \sum_{j=0}^{m-r} 
\frac{B_{j}(\va_{v})}{j!(m-r - j)!}    
\sum_{q \in  \Z^d - tv \; (  \hspace{-.1cm}   \bmod \Pi_v)}       
{\langle  q, z \rangle}^{m-r-j}\\
&=
\sum_{v \in \rm{V}}    
   \langle v, z \rangle^{m-1}       
 \sum_{j=0}^{1} 
\frac{B_{j}(\va_{v})}{j!}            
\sum_{q \in  \Z^d - tv \; (  \hspace{-.1cm}   \bmod \Pi_v)}      
\frac{ {\langle  q, z \rangle}^{1 - j} }{(1 - j)!} \\  \label{step 2 for smooth}
&= \sum_{v \in \rm{V}}       
\frac{    \langle v, z \rangle^{m-1}  }{ \prod_{k=1}^d  \langle w_k(v), z \rangle}         
\sum_{q \in  \Z^d - tv \; (  \hspace{-.1cm}  \bmod \Pi_v)} 
\langle  q, z \rangle        
 - \frac{1}{2} \sum_{v \in \rm{V}} 
 \langle v, z \rangle^{m-1} 
 \vol \Pi_{v}  \frac{   \sum_{k=1}^d \langle w_k(v), z \rangle }{\prod_{k=1}^d \langle w_k(v), z \rangle}.
\end{align}  


In other words, we have
\begin{align} 
\label{vanishing identity for m=d, r=d-1}
 \sum_{v \in \rm{V}}   
 \frac{ \langle v, z \rangle^{m-1}  }{ \prod_{k=1}^d  \langle w_k(v), z \rangle}         
\sum_{q \in  \Z^d - tv \; (  \hspace{-.1cm}   \bmod \Pi_v)} 
\langle  q, z \rangle        
 = \frac{1}{2} \sum_{v \in \rm{V}}   
 \langle v, z \rangle^{m-1} \vol \Pi_{v}        
 \frac{   \sum_{k=1}^d \langle w_k(v), z \rangle }{\prod_{k=1}^d \langle w_k(v), z \rangle},
\end{align}  
for all simple rational polytopes, and for almost all $z\in \C^d$.

Next, we make a simple observation that holds for the family of smooth (integer) polytopes.  If $\P$ is a smooth polytope, then the latter identity
\eqref{vanishing identity for m=d, r=d-1}
enjoys a simplification.  For a smooth polytope, we have
$\sum_{p \in \Pi_v \cap \Z^d}    \langle  p, z \rangle=0$, because the origin is the only integer point in $\Pi_v$. Therefore 
\begin{equation}
\label{Brion-Vergne vanishing}
 \sum_{v \in \rm{V}} 
 \langle v, z \rangle^{m-1}
   \frac{   \sum_{k=1}^d \langle w_k(v), z \rangle }{\prod_{k=1}^d \langle w_k(v), z \rangle}=0,
\end{equation}
for each fixed $1\leq m \leq d-1$, and for almost all $z\in \C^d$.

\begin{remark}
A special case of Corollary \ref{thm:Ehrhart.from.Bernoulli.Barnes} is the case $r=d$. 
 An elementary fact from Ehrhart theory \cite{Sinai2} is that the leading coefficient of $L_\P(t)$ is always equal to $\vol \P$. Recalling that $B_{0}(\va)  = \frac{1}{a_1 \cdots a_d}$, Corollary 
 \ref{thm:Ehrhart.from.Bernoulli.Barnes} tells us that 
\begin{equation}  \label{Brion volume}
\vol \P = \frac{  (-1)^d }{d!}     \sum_{v \in \rm{V}}
   \left(\vol \Pi_{v}\right)        \langle v, z \rangle^{d}  B_{0}(\va_{v})
   =   \frac{  (-1)^d }{d!}     \sum_{v \in \rm{V}}
   \left(\vol \Pi_{v}\right)   \frac{   \langle v, z \rangle^{d} }{  \prod_{m=1}^d \langle w_m(v), z \rangle},
\end{equation}
for almost all $z \in \C^2$.  We've 
recovered here a known volume identity of Brion \cite[Corollary 2]{Brion1} and Lawrence \cite{Lawrence}.  For a Fourier-analytic proof see \cite[Theorem 7.23]{Sinai2}.
 \hfill $\hexago$
\end{remark}





\section{Further remarks and research directions}
\label{sec: further remarks}

We note that for fixed dimension $d$, direct evaluation of the formula for $L_P(t)$ in
Corollary~\ref{thm:Ehrhart.from.Bernoulli.Barnes} need not be a polynomial-time procedure,
because the finite moment sums may contain exponentially many terms as a function of the
binary input of the edge vectors.  
 
\begin{question}
\label{conjecture for computing discrete moments}
Fix the dimension $d$,   and fix a half-open integer parallelepiped $\Pi\subset \R^d$, defined by its edges $w_1, \dots, w_d$.  Fix also a positive integer $k$ and $z \in \C^d$. Let $N$ be the bit input of the edges $w_1, \dots, w_d$.   
Is there an $N^{O(d)}$-time algorithm for computing a polytope Dedekind sum?
\end{question}
Barvinok's algorithm \cite{Barvinok1} guarantees a complexity bound of 
$N^{O(d \log d)}$ for general rational polytopes, where $N$ is the bit input of $\P$.  So in particular Barvinok's algorithm guarantees a polynomial-time algorithm for computing  $\mu_k(\Pi)$ in fixed dimension.  But with half-open integer parallelepipeds we have a lot more structure, suggesting a slightly better bound.

\begin{remark}[A conditional vertex-reduction idea]
\label{remark: variety idea}
Here is a possible way to reduce, for a fixed dilation $t>0$, the number of discrete moments that must be evaluated in Theorem~\ref{thm:Ehrhart quasi-polynomial from Barnes}.  Suppose that $v\in\P$ is a particularly difficult vertex, in the sense that its vertex parallelepiped $\Pi_v$ contains $2^{100}$ lattice points, for example.

Let $d\geq2$, fix the vertex $v$ and the dilation $t>0$, and set
\begin{equation}
D_v(z)
:=
\prod_{j=1}^d\langle w_j(v),z\rangle.
\label{eq:denominator-polynomial-remark}
\end{equation}
By Appendix~\ref{Basic properties of Barnes}, for each $k\geq1$ the product
$
D_v(z)\,
B_k\!\left(t\langle v,z\rangle,\va_v\right)
$
is a polynomial in $z_1,\ldots,z_d$.  We therefore define the algebraic set
\begin{equation}
X_{v,t}
:=
\bigcap_{k=1}^{d-1}
\left\{
z\in\C^d
\;\middle|\;
D_v(z)\,
B_k\!\left(t\langle v,z\rangle,\va_v\right)=0
\right\}.
\label{variety idea}
\end{equation}
Let $\mathcal G_P\subset\C^d$ denote the generic set defined in \eqref{def:almost all z}.  If
$
z\in X_{v,t}\cap\mathcal G_P$,
then $D_v(z)\neq0$, and therefore 
$
B_k\!\left(t\langle v,z\rangle,\va_v\right)=0,
\qquad 1\leq k\leq d-1$.
Consequently, in the vertex-$v$ contribution to Theorem~\ref{thm:Ehrhart quasi-polynomial from Barnes}, part~\eqref{part a of Theorem 2}, all terms involving the intermediate discrete moments
$
\mu_1(\Pi_v,z,tv),\ldots,\mu_{d-1}(\Pi_v,z,tv)
$
vanish. More precisely, the contribution of $v$ reduces to
\begin{equation}
\frac{(-1)^d}{d!}
\left(
B_0(\va_v)\,
\mu_d(\Pi_v,z,tv)
+
B_d\!\left(t\langle v,z\rangle,\va_v\right)
\mu_0(\Pi_v,z,tv)
\right).
\label{eq:conditional-reduced-vertex-contribution}
\end{equation}
Thus, whenever $X_{v,t}$ meets the generic set, we may omit the intermediate moment sums at the chosen vertex when evaluating $L_P(t)$ for that fixed value of $t$.

We note that the conclusion here
 does not eliminate the entire vertex contribution, since the terms containing $\mu_d$ and $\mu_0$ remain.  It also does not eliminate a vertex from the formulas for the Ehrhart coefficients $c_1,\ldots,c_{d-1}$: the set $X_{v,t}$ depends on the fixed dilation $t$, whereas the coefficient formulas involve Barnes numbers and a different collection of moment terms.  The usefulness of \eqref{variety idea} therefore depends on the nonemptiness of $X_{v,t}\cap\mathcal G_P$.
\hfill $\hexago$
\end{remark}

\begin{question}
For a fixed vertex $v$ and a fixed dilation $t>0$, when is
$X_{v,t}\cap\mathcal G_P$ nonempty?  When it is nonempty, can a point in this intersection be found effectively?
\end{question}

\begin{question}
For a fixed dilation $t>0$, can one choose a single generic vector $z$ lying simultaneously in
\[
X_{v_1,t}\cap\cdots\cap X_{v_s,t}
\]
for two or more selected vertices $v_1,\ldots,v_s$, thereby removing all intermediate moment contributions at those vertices?
\end{question}

\begin{question}
Are there additional algebraic conditions on $z$ that can also eliminate the remaining $\mu_d$-term in \eqref{eq:conditional-reduced-vertex-contribution}, and hence remove the entire contribution of a selected vertex for a fixed dilation $t$?
\end{question}

\begin{question}
 Can we understand complete period collapse for $d$-dimensional rational polytopes by understanding the interaction between the lattice flows at each of the vertex parallelepipeds,    appearing in Theorem \ref{thm:Ehrhart quasi-polynomial from Barnes} ?
What is the smallest rational period, over $0<t\leq 1$, of just a single polytope Dedekind sum?
\end{question}

\bigskip
\begin{appendices}

\section{Polytope Dedekind sum reciprocity gives the classical Dedekind sum reciprocity}
\label{sec:dedekind-reciprocity-fractional-part-moments}

Here we show how to retrieve the classical Dedekind sums from the reciprocity law for polytope Dedekind sums, namely Corollary~\ref{cor:ConstantCoeff}.
Let $a,b\in\Z_{>0}$ satisfy $\gcd(a,b)=1$.  
We consider the triangle defined by
$T_{a,b}
:=
\operatorname{conv}
\left\{
(0,0),(b,0),(0,a)
\right\}$. 
For  coprime positive integers $h$ and $k$, the classical Dedekind sum is defined by:
$
s(h,k)
:=
\sum_{r=1}^{k-1}
\left(
\left\{\frac{r}{k}\right\}-\frac12
\right)
\left(
\left\{\frac{hr}{k}\right\}-\frac12
\right)$.
We shall recover the classical reciprocity law:
\begin{equation}
s(a,b)+s(b,a)
=
-\frac14
+
\frac1{12}
\left(
\frac ab+
\frac ba+
\frac1{ab}
\right).
\label{eq:dedekind-reciprocity-goal}
\end{equation}
First, we fix a generic vector $z=(X,Y)\in\C^2$.  At each vertex $v$, we write
\begin{equation}
A_v
:=
\langle w_1(v),z\rangle,
\qquad
B_v
:=
\langle w_2(v),z\rangle,
\qquad
\begin{pmatrix}
\alpha_1(n)\\
\alpha_2(n)
\end{pmatrix}
:=
W_v^{-1}n.
\label{eq:dedekind-local-notation}
\end{equation}
Substituting the explicit two-dimensional Barnes numbers into
Corollary~\ref{cor:ConstantCoeff}, and then applying
Lemma~\ref{lemma:discrete moments via fractional parts} at $t=0$, gives us
\begin{equation}
1
=
\sum_{v\in V(T_{a,b})} C_v,
\label{eq:dedekind-Cv-sum}
\end{equation}
where
\begin{equation}
\begin{aligned}
C_v
:=
\frac12
\sum_{n\in\Pi_v\cap\Z^2}
\Bigg[
{}&
\frac{A_v}{B_v}
B_2\!\left(\{\alpha_1(n)\}\right)
+
\frac{B_v}{A_v}
B_2\!\left(\{\alpha_2(n)\}\right)
\\
&
+
2
\left(
\{\alpha_1(n)\}-\frac12
\right)
\left(
\{\alpha_2(n)\}-\frac12
\right)
\Bigg],
\end{aligned}
\label{eq:dedekind-Cv-formula}
\end{equation}
and
$B_2(x)
:=
x^2-x+\frac16$ 
is the ordinary second Bernoulli polynomial.

\noindent
At $v_0=(0,0)$, the primitive edge vectors are $(1,0)$ and $(0,1)$,
and $\Pi_{v_0}\cap\Z^2=\{0\}$.  Hence
\begin{equation}
C_{v_0}
=
\frac{X}{12Y}
+
\frac{Y}{12X}
+
\frac14.
\label{eq:dedekind-origin-contribution}
\end{equation}

\noindent
At $v_x=(b,0)$, the primitive edge vectors are $w_1(v_x)=(-1,0)$ and 
$w_2(v_x)=(-b,a)$, so that:
\begin{equation}
A_{v_x}=-X \text{ and }
B_{v_x}=aY-bX.
\end{equation}
Since
\begin{equation}
W_{v_x}^{-1}n
=
\begin{pmatrix}
-n_1-\dfrac ba n_2\\[2mm]
\dfrac1a n_2
\end{pmatrix},
\label{eq:dedekind-vx-inverse}
\end{equation}
the $a$ lattice points of $\Pi_{v_x}$ are parametrized by
$0\leq r<a$ with
\begin{equation}
\{\alpha_1(n)\}
=
\left\{-\frac{br}{a}\right\},
\qquad
\{\alpha_2(n)\}
=
\frac ra.
\label{eq:dedekind-vx-fractional-coordinates}
\end{equation}
The Bernoulli distribution relation and the definition of the
Dedekind sum give
\begin{equation}
\sum_{r=0}^{a-1}
B_2\!\left(\frac ra\right)
=
\sum_{r=0}^{a-1}
B_2\!\left(\left\{-\frac{br}{a}\right\}\right)
=
\frac1{6a},
\end{equation}
\begin{equation}
\sum_{r=0}^{a-1}
\left(
\left\{-\frac{br}{a}\right\}-\frac12
\right)
\left(
\frac ra-\frac12
\right)
=
\frac14-s(b,a).
\end{equation}
Therefore,
\begin{equation}
C_{v_x}
=
-\frac{X}{12a(aY-bX)}
-
\frac{aY-bX}{12aX}
+
\frac14
-
s(b,a).
\label{eq:dedekind-vx-contribution}
\end{equation}

\noindent
Finally, at $v_y=(0,a)$, we have
\begin{equation}
w_1(v_y)=(0,-1),
\qquad
w_2(v_y)=(b,-a),
\qquad
(A_{v_y},B_{v_y})=(-Y,bX-aY).
\label{eq:dedekind-vy-data}
\end{equation}
The identical calculation, with $a,b$ and $X,Y$ interchanged, gives
\begin{equation}
C_{v_y}
=
-\frac{Y}{12b(bX-aY)}
-
\frac{bX-aY}{12bY}
+
\frac14
-
s(a,b).
\label{eq:dedekind-vy-contribution}
\end{equation}

Substituting
\eqref{eq:dedekind-origin-contribution},
\eqref{eq:dedekind-vx-contribution}, and
\eqref{eq:dedekind-vy-contribution} into
\eqref{eq:dedekind-Cv-sum}, and simplifying the rational terms, we see that the $z:=(X,Y)$ parameter magically cancels out (as it must).  We therefore obtain the reciprocity law
\eqref{eq:dedekind-reciprocity-goal}:
\begin{equation}
1
=
\frac34
+
\frac{a^2+b^2+1}{12ab}
-
s(a,b)
-
s(b,a).
\label{eq:dedekind-final-vertex-sum}
\end{equation}

\bigskip
\section{Known properties of Barnes polynomials}
\label{Basic properties of Barnes}

Throughout this appendix, let $d\geq 1$ be an integer and let
$\va=(a_1,\ldots,a_d)\in(\C^{*})^d$.
 Unless otherwise stated, the index
$k$ is a nonnegative integer and the scalar variables $t,x,y$ belong to
$\C$. Thus, although the dilation parameter in our Ehrhart-theoretic
applications satisfies $t>0$, no positivity restriction on $t$ is needed
for the polynomial identities in this appendix.

Given their defining generating function
\eqref{Barnes}, it is natural that the Barnes polynomials should inherit some of the properties of the classical Bernoulli polynomials.  First, the Barnes polynomials can be easily
written as a finite sum of classical Bernoulli polynomials: 
\begin{equation} \label{explicit sum}
\frac{1}{k!} B_k(t, \va) = \sum_{m_1 +\cdots+m_d=k} 
   \frac{a_1^{m_1 - 1}}{m_1 !}  \cdots    \frac{a_d^{m_d-1}}{m_d !}
B_{m_1}\Big(\tfrac{t}{a_1 +\cdots+a_d}\Big) \cdots   B_{m_d}\Big(\tfrac{t}{a_1 +\cdots+a_d}   \Big),
\end{equation}
where the $B_m(x)$ are the classical Bernoulli polynomials and the
indices are nonnegative integers $m_j \geq 0$ (see Bayad and Kim
\cite[Theorem 2]{BayadKim2013}).
We see that if all of the $a_j$'s are rational numbers, then $B_k(t, \va)$ is a polynomial in $t$, with rational coefficients.

It follows immediately from 
\eqref{explicit sum}
that the Barnes numbers 
$B_k(\va):= B_k(0, \va)$
may be written in terms of the classical Bernoulli numbers as follows: 
 \begin{equation} \label{Barnes numbers}
B_k(\va) = \sum_{m_1 + \cdots +  m_d=k} 
\binom{k}{m_1, \dots, m_d}
  a_1^{m_1-1}    \cdots  a_d^{m_d-1}
B_{m_1} \cdots   B_{m_d},
\end{equation}
 where the $B_j$'s are the classical Bernoulli numbers 
 \cite[Theorem 2]{BayadKim2013}.  Equation \eqref{Barnes numbers} also shows that $(a_1 \cdots a_d) B_k(\va)$ is a homogeneous polynomial of the $a_j$'s, of degree $k$. 
The Barnes numbers extend the usual Bernoulli numbers, in the precise sense that for dimension $d=1$, and for
$ \va = (1)$, we have $\ B_j(\va) = B_j$, the $j$'th Bernoulli number.
Many of our results used the following explicit polynomial representation for the 
Barnes polynomials in terms of the Barnes numbers 
$B_k(\va)$:
\begin{equation}   \label{explicit polynomial in terms of Barnes numbers}
B_k(t, \va) = \sum_{j=0}^k  
\binom{k}{j}
B_{j}(\va) t^{k-j},
\end{equation}
for each $k \geq 0$.
From  \eqref{Barnes numbers}  and \eqref{explicit polynomial in terms of Barnes numbers}, it follows at once that 
$(a_1 a_2 \cdots a_d) B_k(t, \va) \in
\Q[t,a_1, \dots, a_d]$.
Also, it's clear that the coefficients of $B_k(t, \va)$, as a polynomial in $t$,  are symmetric functions of 
$a_1, \dots, a_d$.  
The Barnes polynomials enjoy the following strong properties. In
parts (a)--(d), let $t,x,y\in\C$ and let
$\lambda\in\C^{\times}$.
\begin{enumerate}[(a)]
\item 

[{\bf Reciprocity}]
\label{reciprocity}
\begin{align*}  
(-1)^k  B_k(-t, \va) &= B_k( t + a_1 + \cdots + a_d, \va). 
\end{align*}
\item 

[{\bf Homogeneity}]
\label{homogeneity}
\begin{equation*}   
B_k(\lambda t, \lambda \va) = \lambda^{k-d}  B_k(t, \va).
\end{equation*}

\item

[{\bf Differentiation}]
\label{differentiation}
For $k\geq 1$,
\begin{equation*}   
\frac{d}{dt} B_k(t, \mathbf{a}) 
= k B_{k-1}(t, \mathbf{a}).
\end{equation*}

\item 

[{\bf Addition formula}]
\label{Addition formula}
\begin{equation*}
  B_k(x + y, \mathbf{a}) 
  = \sum_{i=0}^{k} \binom{k}{i} B_i(x, \mathbf{a}) y^{k-i}.  
\end{equation*} 
\end{enumerate}

\noindent
In the homogeneity formula,
$\lambda\va:=(\lambda a_1,\ldots,\lambda a_d)$.
Reciprocity, homogeneity, and the addition formula hold for every
$k\geq0$, while differentiation holds for every $k\geq1$.
See \cite{BayadKim2013} for  proofs and more background on Barnes polynomials.
For the purposes
 of explicitly computing some Ehrhart polynomials, it is very useful to list the first few Barnes polynomials in low dimensions, relevant for us in the practical computations of the theorems above.  It turns out that  the standard power sums
 \begin{equation}
 s_r:= \displaystyle\sum_{j=1}^d a_j^r
 \end{equation}
 help us write the formulas more compactly, especially as the dimension increases.

\smallskip
\noindent
{\bf The first 7 Barnes polynomials in dimension d}:
\begin{align}
\label{B_0(t, vector a}
 B_0(t, \va) &= \frac{1}{a_1  \cdots a_d}. \\
  \label{B_1(t, vector a}
(a_1  \cdots a_d)
 B_1(t, \va) 
&= 
 t -\frac{1}{2}   
\left(a_1 + a_2 + \cdots + a_d\right)
= t -\frac{1}{2}  s_1.\\  
 \left(a_1  \cdots a_d \right)
 B_2(t, \va) 
&=
\left( t^2 - t \sum_{i=1}^d a_i 
- \frac{1}{12} \sum_{i=1}^d a_i^2 + \frac{1}{4} \left( \sum_{i=1}^d a_i \right)^2 \right)
\\
&=
\left(t-\frac{1}{2}s_1\right)^2
-\frac{1}{12}s_2.
\end{align}

\begin{align}
 (a_1  \cdots a_d )B_3(t, \va) 
&=
 t^{3} - \frac{3}{2}s_{1} t^{2} 
 + \frac{1}{4}\left(3s_{1}^{2} - s_{2}\right)t 
 + \frac{1}{8}\left(s_{1}s_{2} - s_{1}^{3}\right)\\
 &=
\left(t - \frac{1}{2}s_{1}\right)^{3}
-\frac{1}{4} s_{2}\left(t - \frac{1}{2}s_{1}\right).
\end{align}

\begin{align}
(a_1  \cdots a_d )B_4(t, \va)
&=
\left(t - \frac{1}{2}s_{1}\right)^{4} - \frac{1}{2}s_{2} \left(t - \frac{1}{2}s_{1}\right)^{2} + \frac{1}{48}s_{2}^{2} + \frac{1}{120}s_{4}.
\end{align}

\begin{remark}
\label{remark: parity of Barnes polynomials}
It becomes apparent, as suggested by the latter  evaluations of the Barnes polynomials, that the Barnes polynomials enjoy the following property.   
For odd $k$  the Barnes polynomials
$B_k(t, \va)$ are odd polynomials in the variable 
$y:=\left(t-\frac{1}{2} s_1\right)$, while for even $k$ the Barnes polynomials 
$B_k(t, \va)$ are even polynomials in the same variable 
$y:=\left(t-\frac{1}{2} s_1\right)$.  This is very easy to prove, via a symmetry argument in their generating function.  
\end{remark}

\begin{align}
(a_1  \cdots a_d ) B_{5}(t,a)
=
y^{5}
-\frac{5}{6}s_{2} \,y^{3}
+\left(
\frac{5}{48}s_{2}^{2}
+\frac{1}{24}s_{4}
\right)
y.
\end{align}

\begin{align}
(a_1  \cdots a_d ) 
B_6(t,a)
=
y^6
-\frac54s_2 \, y^4 
+
\left(
\frac{5}{16}s_2^2+\frac18s_4
\right)
y^2
-
\left(
\frac{5}{192}s_2^3
+
\frac{1}{32}s_2s_4
+
\frac{1}{252}s_6
\right).
\end{align}


\bigskip
\noindent{\bf The first five Barnes numbers in dimension $d$}:
\begin{align}
\label{B_0(vector a)}
 B_0(\va) 
&= 
\frac{1}{a_1  \cdots a_d}. \\
  \label{B_1(vector a)}
 B_1(\va) 
&=
-\frac{1}{2}   \frac{a_1 + a_2 + \cdots + a_d}{a_1  \cdots a_d}.\\
  \label{B_2(vector a)}
 (a_1  \cdots a_d ) B_2(\va) 
&=
\frac{1}{12} \left(  3 s_1^2-s_2 \right).  \\ \label{B_3(a), dim d}
(a_1  \cdots a_d )B_3(\va) 
&=
\frac{1}{8}(s_{1}s_{2} - s_{1}^{3}). \\
(a_1  \cdots a_d ) B_4(\va)
&=\frac{1}{16} s_{1}^{4} - \frac{1}{8}s_{1}^{2}s_{2} 
+ \frac{1}{48} s_{2}^{2} + \frac{1}{120} s_{4}.
\end{align}


\bigskip
\noindent
For the purposes of computing examples of Ehrhart quasi-polynomials in $\R$ and $\R^2$, we record the explicit special cases of $d=1$ and $d=2$ here.

\bigskip
\noindent
{\bf Barnes polynomials in dimension 1, relevant for us}:
\begin{equation}
  \label{B_1(t, (a_1))}
  B_0(t, (a_1))= \frac{1}{a_1}, \quad  B_1(t, (a_1)) = \frac{1}{a_1} t - \frac{1}{2}. 
\end{equation}


\bigskip
\noindent
{\bf Barnes polynomials in dimension 2, relevant for us}:
\begin{align}
 B_0(t, \va) 
&= 
\frac{1}{a_1 a_2}  \\
 B_1(t, \va) 
&= 
\frac{1}{a_1 a_2} t -\frac{1}{2}   
\frac{a_1 + a_2}{a_1  a_2}.  \\
 B_2\left(t, \va\right)
&= 
t^2 \left(   \frac{1}{a_1 a_2}  \right)
- t \left(   \frac{1}{a_2}   +  \frac{1}{a_1} \right) 
+\frac{1}{6}\left(\frac{a_2}{a_1} + \frac{a_1}{a_2}\right) + \frac{1}{2}. 
\label{B_2(t, (a1, a2))}
\end{align}

\noindent
{\bf Barnes numbers in dimension 2, relevant for us}:
\begin{align}
B_0(\va) 
 = 
\frac{1}{a_1 a_2}, \quad 
B_1(\va) 
= 
 -\frac{1}{2}   
\frac{a_1 + a_2}{a_1  a_2}, \quad   
 B_2(\va)
=   
\frac{1}{6}\left(\frac{a_2}{a_1} + \frac{a_1}{a_2}\right) + \frac{1}{2}. 
\label{B_2((a1, a2))}
\end{align}

\end{appendices}


\bigskip
\noindent
{\bf AI Usage}. \ The author did not use AI tools in the discovery, development, or writing of any of the results or proofs in this paper.  The author did use AI tools to proofread the manuscript and assist in the construction of Example \ref{Example:rational triangle}, which was subsequently checked and verified by the author.

\end{document}